\documentclass[10pt]{amsart}

\usepackage{amsthm,amsmath,amssymb,amscd}
\usepackage{amssymb}
\usepackage{graphics}

\newtheorem{teor}{Theorem}[section]

\newtheorem{defin}{Definition}[section]
\newtheorem{remar}{Remark}[section]
\newtheorem{prop}{Proposition}[section]
\newtheorem{corol}{Corollary}[section]
\newtheorem{lemma}{Lemma}[section]

\newcommand{\intero}{\mathbb{Z}}

\newcommand{\K}{K\"{a}hler}

\newcommand{\R}{\mathbb{R}}

\newcommand{\C}{\mathbb{C}}

\newcommand{\e}{\varepsilon}
\newcommand{\del}{\partial}

\begin{document}

\title[Blowing up and desingularizing K\"{a}hler manifolds]
{Blowing up and desingularizing constant scalar curvature K\"{a}hler
manifolds}
\author[Claudio Arezzo] {Claudio Arezzo}
\address{ claudio.arezzo@unipr.it \\ Universita' di Parma\\Italy}

\author[Frank Pacard] {Frank Pacard}
\address{pacard@univ-paris12.fr\\ University Paris 12 and
Institut Universitaire de France, France}

\maketitle

\vspace{-,15in}


\section{Introduction}

\begin{abstract}
This paper is concerned with the existence of constant scalar
curvature \K\ metrics on blow ups at finitely many points of compact
manifolds which already carry constant scalar curvature \K\ metrics.
We also consider the desingularization of isolated quotient
singularities of compact orbifolds which already carry constant
scalar curvature  \K\ metrics.
\end{abstract}

\medskip

Let $(M,\omega)$ be either a $m$-dimensional compact \K\ manifold or
a $m$-dimensional compact \K\ orbifold with isolated singularities.
By definition, any point $p \in M$ has a neighborhood biholomorphic
to a neighborhood of the origin in ${\mathbb C}^{m} / \, \Gamma$,
where $\Gamma$ is a finite subgroup of $U(m)$ (this last fact is a
consequence of the \K\ property) acting freely on ${\mathbb C}^{m} -
\{ 0 \}$. Observe that, when $p$ is a smooth point of $M$, the group
$\Gamma$ reduces to the identity. In the case where $M$ is an
orbifold, the \K\ form $\omega$ lifts, near any of the singularities
of $M$, to a \K\ form $\tilde \omega$ on a punctured neighborhood of
$0$ in ${\mathbb C}^{m}$. We will always assume that $\tilde \omega$
can be smoothly extended through the origin, i.e. that $\omega$ is
an orbifold metric.

\medskip

If we further assume that the \K\ form $\omega$ has constant scalar
curvature and if we are given $n$ distinct (smooth) points $p_1,
\ldots, p_n \in M$, one of the questions we would like to address in
this paper is whether the blow up of $M$ at the points $p_1, \ldots,
p_n$ can still be endowed with a constant scalar curvature \K\ form.
In this direction, we have obtained the following positive answer~:
\begin{teor}
Let $(M, \omega)$ be a constant scalar curvature compact \K\
manifold or \K\ orbifold with isolated singularities. Assume that
there is no nonzero holomorphic vector field vanishing somewhere on
$M$. Then, given finitely many smooth points $p_1, \ldots, p_n \in
M$ and positive numbers $a_1, \ldots, a_n >0$, there exists $\e_0
>0$ such that the blow up of $M$ at $p_1, \ldots, p_n$
carries constant scalar curvature \K\ forms
\[
\omega_\e \in \pi^{*} \, [\omega] - \e ^{2} \, (a_1 \, PD[E_1] +
\ldots + a_n \, PD[E_n]),
\]
where the $PD[E_j]$ are the Poincar\'e dual of the $(2m-2)$-homology
classes of the exceptional divisors of the blow up at $p_j$ and $\e \in
(0,\e_0)$.

\medskip

If the scalar curvature of $\omega$ is not zero then the scalar
curvatures of $\omega_\e$ and of $\omega$ have the same signs.
\label{th:1.1}
\end{teor}

Following a suggestion of C. LeBrun, we also show that the proof of
Theorem~\ref{th:1.1} can be used to produce zero scalar curvature
\K\ metrics provided the \K\ form $\omega$ we start with has zero
scalar curvature and the first Chern class of $M$ is not zero.
\begin{corol}
\label{co:1.1} Let $(M, \omega)$ be a zero scalar curvature compact
\K\ manifold or orbifold with isolated singularities. Assume that
there is no nonzero holomorphic vector field  vanishing somewhere on
$M$ and that the first Chern class of $M$ is non zero. Then the blow
up of $M$ at finitely many (smooth) points carries zero scalar
curvature \K\ forms.
\end{corol}

Observe that, on manifolds (or orbifolds with isolated
singularities) with discrete automorphism group, there are no
nontrivial holomorphic vector fields. Hence, if these carry a
constant scalar curvature \K\ form, they are examples to whom our
results do apply.

\medskip

On the other hand the assumption is verified also by some manifold
with a continuous family of automorphisms. For example \K\ flat tori
can be used as base manifolds in Theorem~\ref{th:1.1} (but not in
Corollary~\ref{co:1.1} since their first Chern class vanish).

\medskip

Theorem~\ref{th:1.1} and Corollary~\ref{co:1.1} are consequences of
a more general construction which also allows one to desingularize
isolated singularities of orbifolds. This desingularization
procedure combined with Theorem~\ref{th:1.1} is enough to prove the
following~:
\begin{teor}
\label{th:1.2} Any compact complex surface of general type admits
constant scalar curvature \K\ forms.
\end{teor}

It is worth pointing out that some assumption on the initial
manifold $(M, \omega)$ is indeed necessary for either the
desingularization or the blow up procedure to be successful. In
first place we know from the work of Matsushima \cite{matsu} and
Lichnerowicz \cite{li} that the automorphism group of a manifold
with a \K\ constant scalar curvature  metric must be reductive,
hence, for example, the projective plane blown up at one or two
points does not to admit any constant scalar curvature \K\ metric
(see e.g. \cite{besse} page 331). In the same spirit let us recall
that, given a compact complex orbifold $M$ and a fixed \K\ class
$[\omega]$, there is another obstruction for the existence of a
constant scalar curvature \K\ metric in the class $[\omega]$. This
obstruction was discovered by Futaki in the eighties \cite{fu1},
\cite{fu2}, \cite{fu3} for smooth metrics and was extended to
singular varieties by Ding and Tian \cite{dt} and to \K\ constant
scalar curvature metrics by Bando, Calabi \cite{ca2} and Futaki.
This obstruction will be briefly described in Section 4, since it
will play some r\^ole in our construction. The nature of this
obstruction (being a character of the Lie algebra of the
automorphism group) singles out two different types of \K\ manifolds
or \K\ orbifolds with isolated singularities where to look for
constant scalar curvature metrics~: those with no nonzero holomorphic vector fields vanishing somewhere, where the obstruction is vacuous since they do not have any
nontrivial holomorphic vector field, and the others where the Futaki
invariant has to vanish for all holomorphic vector fields.

\medskip

We will say that $(M, \omega)$, a constant scalar curvature \K\
manifold or \K\ orbifold with isolated singularities,  is {\em
nondegenerate} if it does not carry any nontrivial holomorphic
vector field vanishing somewhere (note that this definition does not
depend on the particular \K\ class on $M$), and we will say that
$(M, \omega)$ is {\em Futaki nondegenerate} if the differential of
the Futaki invariant satisfies some nondegeneracy condition (this
definition, which will be made precise in Section 4, does depend on
the \K\ class $[\omega]$).

\medskip

Note that, thanks to Matsushima-Lichnerwicz's decomposition of the
Lie algebra of holomorphic vector fields on a manifold which admits
a constant scalar curvature \K\ metric (see e.g. \cite{besse} and
\cite{fu3}), a  constant scalar curvature \K\ manifold is {\em
nondegenerate} if and only if every holomorphic vector field is
parallel.

\medskip

Our construction gives a quite precise description of the \K\ forms
we obtain on the blown up manifold or on the desingularized
orbifold. We shall now describe more carefully the general
construction and some of its consequences, but also we shall give
more details about the \K\ forms we construct.

\medskip

The construction is obtained by choosing finitely many points $p_1,
\ldots, p_n \in M$ and replacing a small neighborhood of each point
$p_j$, biholomorphic to a neighborhood of the origin in ${\mathbb
C}^{m} / \, \Gamma_j$, by a (suitably scaled down by a small factor
$\e$) piece of a \K\ manifold or a \K\ orbifold with isolated
singularities $(N_j,\eta_j)$, biholomorphic to ${\mathbb C}^{m} / \,
\Gamma_j$ away from a compact subset. This generalized connected sum
yields a \K\ manifold or a \K\ orbifold with isolated singularities
that we call
\[
M \sqcup _{\e , p_1} N_1 \sqcup \ldots \sqcup_{\e,p_n} \, N_n
\]
and whose complex structure does not depend on $\e \neq 0$. We
proceed to perturb the \K\ forms $\omega$ and $\e^{2} \, \eta_j$ on
the various summands, analyzing in Section 5 the linear and in
Section 6 the non linear part of the constant scalar curvature
equation in a given \K\ class. This leads to a study of nonlinear
fourth order elliptic partial differential equations on the \K\
potentials. Then, at the end of Section 6, we ``glue" the \K\
potentials of the perturbed \K\ forms on the different summands to
get a \K\ form whose scalar curvature is constant. The most
important condition that ensures this program to be successful is
the following~: {\em Each $(N_j, \eta_j)$ is an "Asymptotically
Locally Euclidean" (ALE) space and $\eta_j$  is a zero scalar
curvature  \K\ form.}

\medskip

Since the term ALE has often been used with slightly different
meanings, we make precise this definition. In this paper, an ALE
space $(N, \eta)$ is a $m$-dimensional \K\ manifold or \K\ orbifold
with isolated singularities, which is biholomorphic to ${\mathbb
C}^{m} / \, \Gamma$ outside a compact set, where $\Gamma$ is a
finite subgroup of $U(m)$ acting freely on ${\mathbb C}^{m} - \{ 0
\}$, and which is equipped with a \K\ metric $\eta$ which converges
to the Euclidean metric at infinity. In the case where $N$ is an
orbifold with isolated singularities, we assume that, near any
singularity modeled after a neighborhood of $0$ in ${\mathbb C}^{m}/
\tilde \Gamma$, the \K\ form $\eta$ lifts smoothly to a neighborhood
of $0$ in ${\mathbb C}^{m}$. In addition, we will always assume that
there exist complex coordinates $(u^{1}, \ldots, u^{m})$
parameterizing $N$ away from a compact set, in which the \K\ form
$\eta$ can be expanded as
\begin{equation} \eta = i \, \del \, \bar \del \,
( \mbox{$\frac{1}{2}$} \, |u|^{2} + \tilde \varphi(u))
\label{eq:1.1}
\end{equation} at infinity, where the potential $\tilde \varphi$
satisfies
\begin{equation}
\tilde \varphi (u) =  a \, |u|^{4-2m} + {\mathcal O}(|u|^{3-2m}) ,
\label{eq:1.2}
\end{equation}
when $m\geq 3$ and
\begin{equation}
\tilde \varphi (u) =  a  \, \log |u|+ {\mathcal O}(|u|^{-1}) .
\label{eq:1.2bis}
\end{equation}
when $m=2$. Here $a \in {\mathbb R}$ and we agree that ${\mathcal O}
(|u|^{q})$ is a smooth function whose $k$-th partial derivatives are
bounded by a constant times $|u|^{q-k}$, for all $k \geq 2$. The
growth (or decay) of the \K\ potentials for these models is a subtle
problem of independent interest (see e.g. \cite{bkn}). In given
examples, we will see in Section 7 that these potentials can arise
with various orders and decays, and we will show (Lemma
\ref{le:7.2}) that, for zero scalar curvature metrics and under
reasonable growth assumptions on the potential, one can change
suitably the potential in order to get a potential for which
(\ref{eq:1.2})-(\ref{eq:1.2bis}) hold.

\medskip

Let us now summarize the assumptions under which our general
construction works. We will assume that :

\medskip

\begin{itemize}
\item[(i)] $(M,\omega)$ is a $m$-dimensional compact \K\ manifold or
orbifold with isolated singularities.
\\

\item[(ii)] The scalar curvature of $\omega$ is constant.\\

\item[(iii)] $(M, \omega)$ is either {\em nondegenerate} or is {\em Futaki
nondegenerate}. \\

\item[(iv)] Given points $p_1, \ldots, p_n \in M$ which might be either
singular or regular points of $M$, let $\Gamma_j$ be the finite
subgroup of $U(m)$ acting freely on ${\mathbb C}^{m} -\{0\}$ such
that a neighborhood of $p_j$ is biholomorphic to a neighborhood of
the origin in ${\mathbb C}^{m} / \, \Gamma_j$. Each ${\mathbb C}^{m}
/ \, \Gamma_j$ has an ALE resolution $(N_j, \eta_j)$ (which might
either be a manifold or an orbifold with isolated singularities)
endowed with a zero scalar curvature \K\ form $\eta_j$. Furthermore,
we assume that, away from a compact set, the \K\ form $\eta_j$ can
be expanded as in (\ref{eq:1.1}) with a potential satisfying
(\ref{eq:1.2})-(\ref{eq:1.2bis}).
\end{itemize}

\medskip

Our main result reads~:
\begin{teor}
\label{th:1.3} Assume that assumptions (i)-(ii)-(iii) and (iv) are
satisfied. Then, there exists $\e_0 > 0$ and, for all $\e \in (0,
\e_0)$, there exists a constant scalar curvature \K\ form $\tilde
\omega_\e$ defined on $M \sqcup _{{p_{1}, \e}} N_1 \sqcup_{{p_{2},
\e}} \dots \sqcup _{{p_n, \e }} N_n$.

\medskip

As $\e$ tends to $0$, the sequence of \K\ forms $\tilde \omega_\e$
converges (in ${\mathcal C}^{\infty}$ topology) to the \K\ metric
$\omega$, away from the points $p_j$ and the sequence of \K\ forms
$\e^{-2} \, \tilde \omega_\e$ converges (in ${\mathcal C}^{\infty}$
topology) to the \K\ form $\eta_{j}$, on compact subsets of $N_j$.

\medskip

If $\omega$ has positive (resp. negative) scalar curvature then the
\K\ forms $\tilde \omega_\e$ have positive (resp. negative) scalar
curvature.

\medskip

Moreover, if $(M,\omega)$ is {\em nondegenerate}
\[
[\omega_{\e}] = [\omega] + \e^{2}([\eta_1] + \ldots + [\eta_n])
\]
\end{teor}

\noindent Note that, when $(M, \omega)$ is {\em Futaki
nondegenerate}, we cannot control the \K\ class where we find the
constant scalar curvature \K\ metric.

\medskip

All the previous results are consequences of this Theorem.

\medskip

\noindent For example the blow up at smooth points is obtained by
our generalized connected sum construction taking $N_j$ to be the
total space of the line bundle ${\mathcal{O}}(-1)$ over ${{\mathbb
P}}^{m-1}$ (in this case $\Gamma_j =\{id\}$). The key property (iv)
asks for an ALE zero scalar curvature metric $\eta_j$ on
${\mathcal{O}}(-1)$ such that $[\eta_j] = -PD[E_j]$ and with
appropriate decay at infinity. These \K\ forms have been obtained by
Calabi \cite{ca}. When $m=2$, $\eta$ is usually referred to in the
literature as Burns metric and it has been described (and
generalized) in a very detailed way by LeBrun \cite{lb}. In higher
dimensions, $m\geq 3$, these metrics have been generalized by
Simanca \cite{si}. The ALE property and the issue of the rate of
decay of these metrics towards the Euclidean metric can be easily
derived from these papers. In the $2$ dimensional case, the \K\ form
$\eta$ is explicit and these properties follow at once, while, in
higher dimensions, it can be shown that these metrics have a
potential for which (\ref{eq:1.1}) and
(\ref{eq:1.2})-(\ref{eq:1.2bis}) are satisfied. The analysis of
these asymptotic properties will be done in Lemma~\ref{le:7.1} (Raza
\cite{ra} has given an alternative proof using toric geometry). In
any case, assumption (iv) is fulfilled and, given smooth points
$p_1, \dots p_n \in M$ and positive constants $a_1, \ldots, a_n$,
the existence of such models can be plugged into
Theorem~\ref{th:1.3} with all ALE spaces equal to $N =
{\mathcal{O}}(-1)$ over ${{\mathbb P}}^{m-1}$ with the
Burns-Calabi-Simanca form $\eta_j  = a_j \, \eta$. This leads to the
results of Theorem~\ref{th:1.1}, which then also holds for {\em
Futaki nondegenerate} manifolds $(M,\omega)$, only losing control on
the \K\ class to represent.

\medskip

In Section ~8 we will observe that our gluing procedure decreases
the starting scalar curvature. Therefore if $(M,\omega)$ has zero
scalar curvature Theorem~\ref{th:1.3} gives (small) negative scalar
curvature metrics. Nonetheless, if the first Chern class is not
zero, Lebrun-Simanca \cite{lbsm} have shown that there exist nearby
\K\ metrics $\omega_+$ and $\omega_{-}$ of (small) positive and
negative constant scalar curvature respectively. We can then apply
Theorem \ref{th:1.1} to $(M,\omega_+)$ and $(M,\omega_{-})$ to get
positive and negative \K\ metrics of constant scalar curvature on
the blow up. We will show in Section ~8 how this implies Corollary
\ref{co:1.1}, a result which also extends to {\em Futaki
nondegenerate} manifolds with nonzero first Chern class. A similar
result had been previously proved in complex dimension $2$ by
Rollin-Singer \cite{rs} who have shown that one can desingularize
compact orbifolds of zero scalar curvature with cyclic orbifold
groups, keeping the scalar curvature zero, by solving on the
desingularization the hermitian anti-selfdual equation.

\medskip

To prove Theorem \ref{th:1.2} we need to apply Theorem \ref{th:1.3}
more than once. The idea, which comes directly form algebraic
geometry, is to associate to $M$ a (possibly) singular complex
surface $\bar{M}$, such that $M$ is obtained form $\bar{M}$ by
desingularizing and blowing up smooth points a finite number of
times. Algebraic geometry (see e.g. \cite{bpvdv}) tells us that if
$M$ is a surface of general type then $\bar{M}$ (which is called the
pluricanonical model of the minimal model of $M$)
\begin{itemize}
\item[(i)]
is again a complex surface \cite{kod},\\
\item[(ii)]
it has only isolated singular points whose local structure groups
$\Gamma_j$ are in $SU(2)$ \cite{br},\\
\item[(iii)]
the first Chern class of $\bar{M}$ is negative, hence it has only a
discrete group of automorphisms (\cite{kobb}, Theorem 2.1 pag. 82), \\
\item[(iv)]
$\bar{M}$ admits a \K -Einstein orbifold metric \cite{k}.\\
\end{itemize}
We will explain below how Theorem \ref{th:1.3} can be used to
resolve $SU(2)$ singularities. Granted this, $M$ is then reobtained
form this desingularized manifold after a finite number of blow ups
at smooth points, and the constant scalar curvature \K\ metric is
then given by Theorem \ref{th:1.1}.

\medskip

If  $p_j$ is a singular point (and hence $\Gamma_j$ is not the
identity group), there is no unique way to resolve the singularity,
and in fact this is an extremely rich area of algebraic geometry.
Once again, whether constant scalar curvature metrics exist or not
on such resolutions depends, according to Theorem~\ref{th:1.3}, on
the existence of ALE scalar flat \K\  resolution of ${\mathbb C}^{m}
/ \, \Gamma$. For a general finite subgroup $\Gamma \subset U(m)$,
the existence of such a resolution is unknown and this prevents us
to state general existence results for constant scalar curvature \K\
metrics. Nonetheless there are large classes of discrete nontrivial
groups for which a good local model is known to exist, looking at
Ricci-flat metrics, very much in the spirit of non compact versions
of the Calabi conjecture.

\medskip

This is the line started by Tian-Yau \cite{ty1} and Bando-Kobayashi
\cite{bk}, culminating in Joyce's proof of the ALE Calabi conjecture
\cite{j}. Joyce himself used this approach to have good local models
for his well known special holonomy desingularization result.
Joyce's theorem, recalled in Theorem \ref{th:7.1}, implies that
given a ${\mathbb C}^{m} /\Gamma$ such that a \K\ crepant resolution
$N$ exists (i.e. a \K\ resolution with $c_1(N)=0$), then $N$ has a
Ricci-flat \K\ metric $\eta$ which, at infinity, can be expanded as
\[
\eta = i \, \del \, \bar \del \, (\mbox{$\frac{1}{2}$} \, |u|^{2} +
\tilde \varphi(u))
\]
for some potential $\tilde \varphi$ satisfying
\[
\tilde \varphi =  {\mathcal O}(|u|^{2-2m}),
\]
so well inside the range of application of Theorem~\ref{th:1.3}.
This approach works for example when $\Gamma = \intero _m$, acting
diagonally on ${\mathbb C}^{m}$ \cite{ca}, and for {\em any} finite
subgroup of $SU(2)$ and $SU(3)$, since in this case we know that
${\mathbb C}^{m} / \, \Gamma$ has a smooth \K\ crepant resolution
(see \cite{ro} and \cite{j}, Chapter 6.4 and Chapter 8 for these
results). The $2$-dimensional case can be handled directly also
relying on Kronheimer's result \cite{kr}.

\medskip

In light of these results, and when $m = 2$, we can apply
Theorem~\ref{th:1.3} when $(M, \omega)$ is a $2$-dimensional {\em
nondegenerate} or {\em Futaki nondegenerate} \K\ orbifold with
isolated singularities and $p_1, \dots p_n \in M$ is any set of
points with a neighborhood biholomorphic to a neighborhood of the
origin in ${\mathbb C}^{2} / \, \Gamma_j$, where $\Gamma_j$ is a
finite subgroup of $SU(2)$ acting freely on ${\mathbb C}^{2} -
\{0\}$. As explained above, this is enough to prove
Theorem~\ref{th:1.2}.

\medskip

Let us mention that other local models are known, for example when
$\Gamma = {\mathbb Z}_k$ is acting on ${\mathbb C}^{m}$ by
multiplication by a $k$-th root of unity. In this case $N$ is the
total space of the line bundle ${\mathbb O}(-k)$ over ${\mathbb
P}^{m}$ and the metric has zero scalar curvature but, in general, is
not Ricci-flat. Rollin-Singer \cite{rs2} have proved the required
decay properties on the corresponding \K\ potential.

\medskip

Other examples should come from the work of Calderbank-Singer
\cite{cs}. They have in fact shown the existence of ALE zero scalar
curvature resolutions of all $U(2)$ cyclic isolated singularities.
The only piece of information missing at the moment to use them in
our construction is the behavior at infinity of a \K\ potential
associated to these ALE metrics.

\medskip

One of the main source of interest in \K\ metrics with constant
scalar curvature lies in its relation with algebraic geometric
properties of the underlying manifold, such as Chow-Mumford, Tian,
or asymptotic Hilbert-Mumford stability (see for example \cite{t1},
\cite{do}, \cite{pt}, \cite{p}, \cite{ma1}). We know for example, by
Mabuchi's extension of Donaldson's work \cite{ma1}, \cite{do}, that
if an integral class is represented by an extremal \K\ metric, then
the underlying algebraic manifold is asymptotically stable in a
sense which depends on the structure of the automorphism group which
preserves the class. In particular if we have a \K\ manifold with
discrete automorphism group, this stability reduces to the classical
Chow stability. Rescaling the \K\ class $[\omega_{\e}]$ by a factor
$k$ to make it integral, our results have the following corollary~:
\begin{corol}
\label{co:1.2} Let $(M,L)$ be a polarized compact algebraic manifold
of complex dimension $m \geq 2$, with discrete automorphisms group,
and $\omega$ a \K\ form with constant scalar curvature in an
integral class. Then all the manifolds obtained by the blowing up
any finite set of points are asymptotically Chow stable relative to
the polarizing class $k\pi^{*} [\omega]  -  \, (PD[E_1] + \cdots +
\, PD[E_n])$, where the $E_j$ are the exceptional divisors of the
blow up and $k$ is sufficiently big.
\end{corol}

Note that playing with the weights $a_j$ in Theorem \ref{th:1.1} one
gets abundance of different polarizing classes for which the above
corollary holds (in the above statement we have used $a_1 = \cdots =
a_n$). Moreover similar results for the different versions of
stability which are known to be implied by the existence of constant
scalar curvature \K\ metrics follow from our theorems. In this
setting it is interesting to observe that Mabuchi and Donaldson
results do not apply to \K\ orbifolds, due to the failure of
Tian-Catling-Zelditch expansion \cite{so}. Nonetheless if a {\em
full} desingularizing process were to go through, then we would get
the stability of the smooth polarized manifold obtained.

\medskip

Another important phenomenon concerning constant scalar curvature
\K\ metrics is that they are unique in their class, up to
automorphisms. This result was proved for \K -Einstein metrics
thanks to the work of Calabi \cite{cal} when $c_1 \leq 0$ and
Bando-Mabuchi \cite{bm} when $c_1 >0$. Uniqueness of constant scalar
curvature \K\ metrics was then proved by Donaldson and Mabuchi
\cite{ma2} (for extremal metrics) in integral classes (with either
discrete or continuous automorphism group), and by Chen \cite{c} for
any \K\ class on manifolds with $c_1 \leq 0$. Recently  Chen and
Tian \cite{ct} have proved it for any \K\ class even for extremal
metrics. This implies that all the constant scalar curvature \K\
metrics produced in this paper are the unique such representative of
their \K\ class, up to automorphisms.

\medskip

Despite this monumental work, our knowledge of concrete examples is
still limited and mainly confined in complex dimension $2$. For
example, Hong \cite{h1}, \cite{h2} has proved the existence of such
metrics in some \K\ classes of ruled manifolds, and Fine \cite{f}
has studied this problem for complex surfaces projecting over
Riemann surfaces with fibres of genus at least two. The only case
completely understood, and giving a rich source of examples, is the
one of zero scalar curvature \K\ surfaces thanks to the work of Kim,
LeBrun, Pontecorvo and Singer \cite{lbsing}, \cite{lb2},
\cite{klbp}, \cite{kp}, and the recent results of Rollin-Singer
\cite{rs}. Our construction allows one to produce many new constant
scalar curvature \K\ manifolds.

\medskip

In our construction, it is possible to keep track of the geometric
meaning of the parameter $\e$. Indeed, for the blow up construction
$a_j^{m-1}\e^{2m-2}$ gives the volume of the exceptional divisor
$E_j$ (up to a universal constant depending only on the dimension).
The r\^ole of $\e$ in our results has a direct analogue in
Fine-Hong's papers \cite{f}, \cite{h1}, \cite{h2}, replacing the
exceptional divisor with the fiber of the projection onto the
Riemann surface or the projectivized fiber of the vector bundle.

\medskip

The last source of interest in our results we would like to mention
is that they give a reverse picture to Tian-Viaclovsky \cite{tv} and
Anderson \cite{an} study of degenerations of critical metrics, in
the special case of constant scalar curvature \K\ metrics in (real)
dimension $4$. If we look at the sequence of \K\ forms $\omega_\e$
on $M \sqcup _{{p_{1}, \e}} N_1 \sqcup_{{p_{2},\e}} \dots \sqcup
_{{p_n, \e}} N_n$ seen as a fixed smooth manifold, we get from our
analysis that they degenerate in the Gromov-Hausdorff sense to the
orbifold $(M, \omega)$ so they provide many examples of the
phenomenon studied in these works.

\medskip

A natural generalization of these results is to look for gluing
theorems for \K\ metrics with nontrivial automorphisms. The
technical difficulties of this extension give rise to some new
interesting phenomena and will be the subject of a forthcoming paper
\cite{Are-Pac-2}.

\medskip

After the first version of these results were posted in electronic
form, Claude LeBrun has indicated some implications of our main
theorem we first missed (notably Corollary~\ref{co:1.1} and
Theorem~\ref{th:1.2}). We wish to thank him for his suggestions.

\section{Gluing the orbifold and the ALE spaces together}

We start by describing the \K\ orbifold near each of its
singularities and we proceed with a description of the ALE spaces
near infinity.

\medskip

Let $(M , \omega)$ be a  $m$-dimensional \K\ manifold or \K\
orbifold with isolated singularities. We choose points $p_1, \ldots,
p_n \in M$. By assumption, near $p_j$, the orbifold $M$ is
biholomorphic to a neighborhood of $0$ in ${\mathbb C}^{m} /
\Gamma_j$ where $\Gamma_j$ is a finite subgroup of $U(m)$ acting
freely on ${\mathbb C}^{m} -\{0\}$. The group $\Gamma_j$ depends on
the point $p_j$ as the subscript is meant to remind the reader. In
the particular case where $p_j$ is a regular point of $M$, the group
$\Gamma_j$ reduces to the identity.

\medskip

We can choose complex coordinates $z : = (z^{1}, \ldots, z^{m})$ in
a neighborhood of $0$ in ${\mathbb C}^{m}$ to parameterize a
neighborhood of $p_j$ in $M$ and, in these coordinates, the \K\ form
$\omega$ can be expended as
\begin{equation}
\omega = \frac{i}{2} \, \sum_{a} \, dz^{a} \wedge d\bar z^{a} +
\sum_{a,b} \, {\mathcal O}_{j,a,b} (|z|^{2} ) \, dz^{a} \wedge d\bar
z^{b} . \label{eq:2.1}
\end{equation}
near $0 \in {\mathbb C}^{m}$ (see \cite{Gri-Har}). The complex
valued functions ${\mathcal O}_{j,a,b} (|z|^{2})$ are smooth
functions which depend on $j$, $a$ and $b$, vanish at the origin and
whose first order partial derivatives also vanish at the origin.
Even though the coordinates $z$ do depend on $p_j$, we shall not
make this dependence apparent in the notation and we hope that the
meaning will be clear from the context.

\medskip

It will be convenient to denote by
\begin{equation}
\begin{array}{rlll}
B_{j,r} & : = & \{  z \in {\mathbb C}^{m} / \, \Gamma_j
\quad : \quad  |z| < r \}, \\[3mm]
B_{j,r}^{*}  & : = & \{  z \in {\mathbb C}^{m} / \, \Gamma_j \quad
: \quad 0 < |z| < r \} , \\[3mm]
\bar B_{j,r}^{*}  & : = & \{  z \in {\mathbb C}^{m} / \, \Gamma_j
\quad : \quad 0 < |z| \leq  r \} ,
\end{array}
\label{eq:2.2}
\end{equation}
the open ball, the punctured open ball and the punctured closed ball
of radius $r$ centered at $p_j$ (in the above defined coordinates
which parameterize a neighborhood of $p_j$ in $M$). We define, for
all $r >0$ small enough (say $r \in (0, r_0)$)
\begin{equation}
M_r : = M - \cup_j \, B_{j,r}  . \label{eq:2.3}
\end{equation}
In other words, $M_r$ is obtained from  $M$  by excising small ball
centered at the points $p_j$. The boundaries of $M_{r}$ will be
denoted by $\del B_{1, r}, \ldots, \del B_{n,r}$.

\medskip

As promised, we now turn to the description of the ALE spaces near
infinity. We assume that, for each $j = 1, \ldots, n$, we are given
a $m$-dimensional \K\ manifold or \K\ orbifold with isolated
singularities $(N_j, \eta_j)$, with one end biholomorphic to a
neighborhood of infinity in ${\mathbb C}^{m} / \, \Gamma_j$. We
further assume that the \K\ metric $g_j$, which is associated to \K\
form $\eta_j$, converges at order $2-2m$ towards the Euclidean
metric. These assumptions imply that one can choose complex
coordinates $u : = (u^{1}, \ldots, u^{m})$ defined away from a
neighborhood of $0$ in ${\mathbb C}^{m}$ to parameterize a
neighborhood of infinity in $N_j$ and, in these coordinates, the \K\
form $\eta_j$ can be expanded as
\begin{equation}
\eta_j = \frac{i}{2} \, \sum_a \, du^{a} \wedge d\bar u^{a} +
\sum_{a,b} \, {\mathcal O}_{j,a,b} (|u|^{2-2m}) \, du^{a} \wedge
d\bar u^{b} ,
\label{eq:2.2bis}
\end{equation}
away from a  fixed neighborhood of the origin in ${\mathbb C}^{m}$.
Here, the complex valued function ${\mathcal O}_{j ,a,b}
(|u|^{2-2m})$ is a smooth function which depends on $j$, $a$ and
$b$, is bounded by a constant times $|u|^{2-2m}$ and whose $k$-th
partial derivatives  are bounded by a constant (depending on $k$)
times $|u|^{2-2m -k}$. As will be explained in Section 7, this decay
assumption is a natural one and, under some mild assumption, one can
prove that this rate of decay in indeed achieved.

\medskip

It will be convenient to denote by
\begin{equation}
\begin{array}{rllll}
C_{j , R} & : = & \{ u \in {\mathbb C}^{n} / \, \Gamma_j \quad
: \quad |u | >  R \}, \\[3mm]
\bar C_{j , R} & : = & \{ u \in {\mathbb C}^{n} / \, \Gamma_j \quad
: \quad |u | \geq   R \},
\end{array}
\label{eq:2.4}
\end{equation}
the complement of a close large ball and the complement of an open
large ball in $N_j$ (in the  coordinates which parameterize a
neighborhood of infinity in $N_j$). We define, for all $R > 0$ large
enough (say $R
> R_0$)
\begin{equation}
N_{j, R} : = N_j - C_{j , R}. \label{eq:2.5}
\end{equation}
which corresponds to the manifold $N_j$ whose end has been
truncated. The boundary of $N_{j,R}$ is denoted by $\del C_{j,R}$.

\medskip

We are now in a position to describe the generalized connected sum
construction. For all $\e > 0$ small enough, we define a complex
manifold by removing from $M$ small balls centered at the points
$p_j$, for $j= 1, \ldots , n$, and by replacing them by properly
rescaled versions of the ALE space $N_j$. More precisely, for all
$\e \in (0, r_0/R_0)$, we choose $r_\e \in (\e \, R_0 , r_0)$ and
define
\begin{equation}
R_\e : =  \frac{r_\e}{\e} . \label{eq:2.6}
\end{equation}
By construction
\[
M_\e : = M \sqcup _{{p_{1}, \e}} N_1 \sqcup_{{p_{2},\e}} \dots
\sqcup _{{p_n, \e}} N_n ,
\]
is obtained by connecting  $M_{r_\e}$ with the truncated ALE spaces
$N_{1, R_\e}, \ldots, N_{n , R_\e}$. The identification of the
boundary $\del B_{j , r_\e}$ in $M_{r_\e}$ with the boundary $\del
C_{j , R_\e}$ of $N_{j, R_\e}$ is performed using the change of
variables
\[
(z^{1} , \ldots, z^{m} )  = \e \, (u^{1} , \ldots, u^{m}) ,
\]
where $(z^{1}, \ldots, z^{m} )$ are the coordinates in $B_{j , r_0}$
and $(u^{1}, \ldots, u^{m})$ are the coordinates in $C_{j , R_0}$.
Observe that, when all singularities of $M$ are in the set $\{p_1,
\ldots, p_n \}$ and the $N_j$ are all  smooth manifolds, then $M_\e$
is a manifold, otherwise $M_\e$ is still an orbifold.

\section{Weighted spaces}

In this section, we describe weighted spaces on $(M^*, \omega)$
where
\begin{equation}
M^{*} : = M - \{ p_j , \quad : \quad j = 1, \ldots, n \},
\label{eq:3.1}
\end{equation}
as well as weighted spaces on each $(N_j , \eta_j)$.

\medskip

To begin with, we define the weighted space on $(M^{*}, \omega)$.
These weighted spaces are by now well known and have been
extensively used in many connected sum constructions. Roughly
speaking, we are interested in functions whose rate of decay or blow
up near any of the points $p_j$ is controlled by a power of the
distance to $p_j$. To make this definition precise, we first need to
define~:
\begin{defin}
Given $\bar r >0$, $k \in {\mathbb N}$, $\alpha \in (0,1)$ and
$\delta \in {\mathbb R}$, the space ${{\mathcal C}}^{k,
\alpha}_{\delta} (\bar B_{j , \bar r}^{*})$ is defined to be the
space of functions $\varphi \in {\mathcal C}^{k, \alpha}_{loc} (\bar
B_{j , \, \bar r}^{*} )$ for which the norm
\[
\| \varphi \|_{{{\mathcal C}}^{k, \alpha}_{\delta} (\bar B_{j , \bar
r}^{*} )} : = \sup_{0 < r \leq \bar r}  r^{-\delta} \, \| \varphi (r
\, \cdot) \|_{{{\mathcal C}}^{k, \alpha} ( \bar B_{j , 1} -B_{j
,1/2})}
\]
is finite. \label{de:3.1}
\end{defin}
Observe that the function \[ z \longmapsto |z|^{\delta'},
\]
belongs to ${{\mathcal C}}^{k, \alpha}_{\delta} (\bar B_{j , \bar
r}^{*} )$ if and only if $\delta' \leq \delta$. This being
understood, we have the~:
\begin{defin}
Given $k \in {\mathbb N}$, $\alpha \in (0,1)$ and $\delta \in
{\mathbb R}$, the weighted space ${\mathcal C}^{k, \alpha}_{\delta}
(M^{*})$ is defined to be the space of functions $\varphi \in
{{\mathcal C}}^{k, \alpha}_{loc} (M^{*})$ for which the following
norm
\[ \| \varphi \|_{{{\mathcal C}}^{k, \alpha}_{\delta} (M^{*})}  : = \| w
\|_{{{\mathcal C}}^{k, \alpha} (M_{r_0/2})} + \sum_j \, \| \varphi
|_{\bar B_{j ,r_0}^{*}} \|_{{{\mathcal C}}^{k, \alpha}_\delta (\bar
B_{j ,r_0}^{*} )}
\]
is finite. \label{de:3.2}
\end{defin}

With this definition in mind, we can now give a quantitative
statement about the rate of convergence of potential associated to
$\omega$ toward the potential associated to the standard \K\ form on
${\mathbb C}^{m}$, at any of the the points $p_j$. More precisely,
near $p_j$ we can write
\begin{equation}
\omega =  i \, \del \, \bar \del \, (\mbox{$\frac{1}{2}$} \, |z|^{2}
+ \varphi_j ), \label{eq:3.2}
\end{equation}
where $\varphi_j$ is a function which lifts smoothly to a
neighborhood of $0$ in ${\mathbb C}^{m}$.

\medskip

We claim that, without loss of generality, it is possible to choose
the potential $\varphi_j$ in such a way that $\varphi_j \in
{\mathcal C}^{4, \alpha}_4 (\bar B_{j, r_0}^{*} )$ (more precisely,
$\varphi_j \in {\mathcal C}^{4, \alpha} (\bar B_{j, r_0})$ and has
all its partial derivatives up to order $3$ which vanish at $0$).
Indeed, the potential $\varphi_j$ lifts to a smooth potential
defined on a neighborhood of $0$ in ${\mathbb C}^{m}$. We can then
perform the Taylor expansion of this potential at $0$, namely
\[
\varphi_j = \sum_{k=0}^{3}  \varphi_j^{(k)} +  \varphi_j'
\]
where the polynomial $\varphi_j^{(k)}$ is homogeneous of degree $k$
and $\varphi_j'$, together with its partial derivatives up to order
$3$, vanish at $0$. Obviously $\varphi_j^{(0)}$ and
$\varphi_j^{(1)}$ are not relevant for the computation of the \K\
form $\omega$ since $\del \, \bar \del \, (\varphi_j^{(0)}+
\varphi_j^{(1)}) =0$, hence we might as well assume that
$\varphi_j^{(0)} \equiv 0$ and $\varphi_j^{(1)} \equiv 0$. Next,
making use of the fact that the coordinates $(z^{1}, \ldots, z^{m})$
are chosen so that (\ref{eq:2.1}) holds, we see that
\[
\del \, \bar \del \, (\varphi_j^{(2)}+ \varphi_j^{(3)}) = {\mathcal
O} (|z|^{2})
\]
but $\del \, \bar \del \, \varphi_j^{(2)}$ and $\del \, \bar \del \,
\varphi_j^{(3)}$ being homogeneous polynomial of degree $0$ and $1$
respectively, we conclude that $\del \, \bar \del \,
(\varphi_j^{(2)}+ \varphi_j^{(3)}) \equiv 0$. Considering
$\varphi_j'$ instead of $\varphi_j$, we have found a potential which
satisfies the desired property.

\medskip

Similarly, we define weighted spaces on the ALE spaces $(N_j,
\eta_j)$. This time we are interested in functions which decay or
blow up near the end of $N_j$ at a rate which is controlled by a
power of the distance from a fixed point in $N_j$. To be more
specific, we first define~:
\begin{defin}
Given $\bar R >0$, $k \in {\mathbb N}$, $\alpha \in (0,1)$ and
$\delta \in {\mathbb R}$, the space ${{\mathcal C}}^{k,
\alpha}_{\delta} (\bar C_{j ,\bar R})$ is defined to be the space of
functions $\varphi \in {\mathcal C}^{k, \alpha}_{loc}  (\bar C_{j,
\bar R})$ such that the following norm
\[
\| \varphi \|_{{{\mathcal C}}^{k, \alpha}_{\delta} (\bar C_{j , \bar
R})} : = \sup_{R \geq \bar R}  R^{-\delta} \, \| \varphi (R \,
\cdot) \|_{{{\mathcal C}}^{k, \alpha} ( \bar B_{j  , 2} - B_{j ,
1})} ,
\]
is finite. \label{de:3.3}
\end{defin}
Again, observe that the function \[ u \longmapsto |u|^{\delta'}
\]
belongs to ${{\mathcal C}}^{k, \alpha}_{\delta} (\bar C_{j  , \bar
R})$ if and only if $\delta' \leq \delta$. We can now define~:
\begin{defin}
Given $k \in {\mathbb N}$, $\alpha \in (0,1)$ and $\delta \in
{\mathbb R}$, the weighted space ${{\mathcal C}}^{k,
\alpha}_{\delta} (N_j)$ is defined to be the space of functions
$\varphi \in {{\mathcal C}}^{k, \alpha}_{loc} (N_j)$ for which the
following norm
\[
\| \varphi \|_{{{\mathcal C}}^{k, \alpha}_{\delta} (N_j)}  : = \|
\varphi \|_{{{\mathcal C}}^{k, \alpha} (N_{j ,2 R_0})} +  \| \varphi
|_{\bar C_{j  , R_0}} \|_{{{\mathcal C}}^{k, \alpha}_\delta (\bar
C_{j , R_0})}
\]
is finite. \label{de:3.4}
\end{defin}
We can now explain the assumption on the rate of convergence at
infinity of the \K\ form $\eta_j$ toward the standard \K\ form on
${\mathbb C}^{m}$. We will assume that, away from a compact set in
$N_j$
\begin{equation}
\eta_j = i \, \del \, \bar \del \, (\mbox{$\frac{1}{2}$} \, |u|^{2}
+ \tilde \varphi_j) , \label{eq:3.3}
\end{equation}
for some potential $ \tilde \varphi_j$ which satisfies
\begin{equation}
\tilde \varphi_j  -  a_j \, | \cdot |^{4-2m} \in {\mathcal C}^{4,
\alpha}_{3-2m} (\bar C_{j  , R_0}) , \label{eq:3.4}
\end{equation}
when $m\geq 3$ and
\begin{equation}
\tilde \varphi_j (u) + a_j  \, \log | \cdot | \in {\mathcal C}^{4,
\alpha}_{-1} (\bar C_{j  , R_0}), \label{eq:3.4bis}
\end{equation}
when $m=2$, for some $a_j \in {\mathbb R}$. As already mentioned in
the introduction, this is a rather natural assumption which is
fulfilled in many important examples.

\begin{remar}
We  will show in Section 7 that, if one simply assumes that the
potential $\tilde \varphi_j$ associated to $\eta_j$ satisfies
\[
\tilde \varphi_j \in {\mathcal C}^{4, \alpha}_{2-\gamma} (\bar C_{j
, R_0})
\] for some $\gamma >0$, then one can always replace $\tilde \varphi_j$
by some potential $\tilde \varphi_j'$ satisfying
(\ref{eq:3.4})-(\ref{eq:3.4bis}). \label{re:3.1}
\end{remar}

\section{The geometry of the equation}

The material contained in this section is well known (see for
example \cite{fu3}); we include it for completeness and to introduce
the reader to the objects entering into the proofs of our results.
Recall that $(M,\omega)$ is a $m$-dimensional compact \K\ manifold
or a \K\ orbifold with isolated singularities. We will indicate by
$g$ the Riemannian metric associated to $\omega$, $\mbox{Ric}_g$ its
Ricci tensor, $\rho_g$ the Ricci form, and ${\bf s} (\omega)$ its
scalar curvature.

\medskip

Following \cite{ls} and \cite{besse}, we want to understand the
behavior of the scalar curvature under deformations of the \K\ form
$\omega$ of the form
\[
\tilde \omega : = \omega  + i \, \del \, \bar \del \, \varphi +
\beta,
\]
where $\beta$ is a closed $(1,1)$ form and $\varphi$ a function
defined on $M$. In local coordinates $(v^{1}, \ldots, v^{m})$, if we
write
\[
\tilde \omega  =  \frac{i}{2} \, \sum_{a, b} \tilde g_{a \bar b} \,
d v^{a} \wedge d \bar v^{b}
\]
then the scalar curvature of $\tilde \omega$ is given by
\begin{equation}
{\bf s} (\tilde \omega) = - \sum_{a, b} \, \tilde g^{a \bar b} \,
\del_{v^{a}}\, \del_{\bar v^{b}} \, \log \, (\mbox{det} \left(
\tilde g )\right) .\label{eq:4.0} \end{equation} The following
result is proven in \cite{ls}, \cite{besse} Lemma 2.158~:
\begin{prop}
The scalar curvature of $\tilde \omega$ can be expanded in terms of
$\beta$ and $\varphi$ as
\[
{\bf s}  \, (\tilde \omega ) = {\bf s} \,  (\omega) - \left(
\frac{1}{2} \, \Delta^{2}_g \, \varphi + \mbox{Ric}_g \cdot
\nabla^{2}_g \, \varphi + \Delta_g (\omega , \beta) + 2 \, (\rho_g ,
\beta)  \right) + Q_g (\nabla^2 \varphi, \beta),
\]
where $Q_g$ collects all the nonlinear terms and where all operators
on the right hand side of this identity are computed with respect to
the \K\ metric $g$. \label{pr:4.1}
\end{prop}
Being a local calculation this formula holds for orbifolds with
isolated singularities too. Of crucial importance will be the two
linear operators which appear in this formula. First, we set
\begin{equation}
{\mathcal L}_g  : =  \Delta_g (\omega , \cdot ) + 2 \, (\rho_g ,
\cdot ) ,\label{eq:4.1}
\end{equation}
which is a linear operator acting on closed $(1,1)$ forms and we
also define the operator
\begin{equation}
{\mathbb L}_g  : = \frac{1}{2} \, \Delta^{2}_g  + \mbox{Ric}_g \cdot
\nabla_g^{2}, \label{eq:4.2}
\end{equation}
which acts on functions.

\medskip

For a general \K\ metric it can be very difficult to analyze these
operators. Nevertheless geometry comes to the rescue at a constant
scalar curvature metric. Indeed, in this case we have
\begin{equation}
{\mathbb L}_g  = 2 \, ({\bar{\partial}}
\partial^{\#}_g)^{\star} \, ({\bar{\partial}}\partial^{\#}_g ),
\label{eq:4.3}
\end{equation}
where $\partial^{\#}_g \, \varphi $ denotes the $(1,0)$-part of the
$g$-gradient of $\varphi$. In other words,
\[
\del^{\#}_g : =  (\bar \partial \, \cdot \, )^{\#}_g ,
\]
where $\#$ is the inverse of
\[
\begin{array}{rclcllll}
\flat : & TM  \otimes {\mathbb C} & \longrightarrow & T^{*}M \otimes
{\mathbb C} \\[3mm]
& \Xi & \longmapsto & g(\Xi, \, \cdot \, ) .
\end{array}
\]
Using (\ref{eq:4.3}) one observes that to any element $\varphi \in
\mbox{Ker} \, {\mathbb L}_g$ one can associate a holomorphic vector
field, namely $\partial^{\#}_g \varphi $, which vanishes somewhere
on $M$. Indeed, just multiply ${\mathbb L}_g \, \varphi = 0$ by
$\varphi$ and integrate the result over $M$ using (\ref{eq:4.3}) to
conclude that
\[ \int_M |{\bar{\partial}} \,
\partial^{\#}_g \, \varphi|^2 \, dv_g=0 ,
\]
where $dv_g$ denotes the volume form associated to $g$. Therefore,
${\bar{\partial}} \, (\partial^{\#}_g \, \varphi)  =0$. Conversely,
given $\Xi$, a holomorphic vector field vanishing somewhere on $M$,
then $\Xi = J \,v + i \, v$, where $v$ is a Killing vector field,
which also must vanish somewhere. By Proposition 1 in \cite{ls}, $v
= \mbox{Im} \, ( \del^\#_g \, \xi )$ for some real valued function
$\xi$, and hence $\Xi = \del^\#_g \, \xi $. Since $\Xi$ is
holomorphic and $g$ is assumed to have constant scalar curvature, we
conclude that $\xi \in \mbox{Ker} \, {\mathbb L}_g$.

\medskip

For constant scalar curvature \K\ metrics, ${\mathfrak h}_0(M)$ the
space of holomorphic vector fields vanishing somewhere on $M$ is
therefore in one to one correspondence with the nontrivial kernel of
the operator ${\mathbb L}_g$ (i.e. elements of the kernel of
${\mathbb L}_g$ whose mean over $M$ is $0$).

\medskip

We define the nonlinear mapping
\[
\begin{array}{rclcll}
S_{\omega}  : & {\mathcal C}^{4, \alpha} (M) & \longrightarrow
&  {\mathcal C}^{0, \alpha} (M) / \, {\mathbb R} \\[3mm]
& \varphi & \longmapsto & {\bf s} (\omega + i \, \partial
{\bar{\partial}} \, \varphi ) \quad \mbox{mod cte}
\end{array}
\]
This is the set-up for LeBrun-Simanca's implicit function theorem
applied to the mapping $S_\omega$. The application of the implicit
function theorem is based on the result which extends immediately to
orbifolds with isolated singularities.
\begin{prop} \cite{ls}
\label{pr:4.2} Assume that $(M, \omega)$ is {\em nondegenerate} and
further assume that its scalar curvature ${\bf s} (\omega)$ is
constant. Then, the operator
\[
\varphi \longmapsto DS_{\omega}|_{0} \, \varphi = - \, {\mathbb L}_g
\, \varphi
\]
is surjective and has a kernel which is spanned by constant
functions.
\end{prop}

\medskip

Indicating by $\psi_g$ the function (up to constants) which gives
\[
\rho_g = \rho_g^{h} + i \, \partial \, {\bar{\partial}} \, \psi_g
\]
(where $\rho^{h}_g$ is the harmonic representative for $[\rho_g]$)
and by $\Xi \in {\mathfrak h}(M)$ a holomorphic vector field, we can
define the Futaki invariant
\[
{\mathcal F} (\Xi , [\omega]) : = \int _M  \Xi \, \psi _g  \, dv_g,
\]
where $dv_g$ denotes the volume form associated to $g$. Let us
denote by ${\mathfrak h}_0(M)$, the space of holomorphic vector
fields which vanish somewhere in $M$ (which is by the baove discussion a 
linear subspace of ${\mathfrak h}(M)$). By definition, $(M, \omega)$
is {\em Futaki nondegenerate} if the ``linearization" of the Futaki
invariant
\[
D{\mathcal F}_{[\omega]} : {\mathfrak h}_0 (M)
\longrightarrow (H^{(1,1)}(M, {\mathbb C}))^{\star} \] is injective.
It is a standard fact, though non obvious, that ${\mathcal F} (\Xi ,
[\omega])$ only depends on the \K\ class and does not depend on its
representative. On the other hand, if $[\omega]$ has a
representative with constant scalar curvature, then ${\mathcal F}
(\Xi , [\omega])$ vanishes for any $\Xi  \in {\mathfrak h}(M)$.

\medskip

Now define the nonlinear mapping
\[
\begin{array}{rclcll}
\hat S_{\omega}  : &   {\mathcal C}^{4, \alpha} (M) \times
{\mathcal{H}}^{1,1}(M, {\mathbb C}) & \longrightarrow &  {\mathcal
C}^{0, \alpha} (M)/ \, {\R} \\[3mm]
& (\varphi , \beta) & \longmapsto & {\bf s}  (\omega +  i \,
\partial {\bar{\partial}} \, \varphi + \beta ) \quad \mbox{mod cte}
\end{array}
\]
where ${\mathcal{H}}^{1,1}(M, {\mathbb C})$ is the space of
$\omega$-harmonic $(1,1)$ forms. The result of \cite{ls} again
extends to orbifolds with isolated singularities, and we have the~:
\begin{prop}
\cite{ls} Assume that $(M, \omega)$ is {\em Futaki nondegenerate}
and further assume that its scalar curvature ${\bf s} (\omega)$ is
constant. Then, the operator \[ (\varphi, \beta) \longmapsto
D\hat{S}_{\omega} |_{(0,[0])} (\varphi, \beta) = - \, ({\mathbb L}_g
\, \varphi + {\mathcal L}_g \, \beta )
\]
is surjective. \label{pr:4.3}
\end{prop}

\section{Mapping properties}

We construct right inverses for the operator ${\mathbb L}_g$ which
has been defined in the previous section. We first explain this
construction in the case where the bounded kernel of ${\mathbb L}_g$
is spanned by constant functions (this is the case if $(M, \omega)$
is {\em nondegenerate}, i.e. when there are no nontrivial
holomorphic vector field vanishing somewhere on $M$) and then we
will explain the modifications which are needed to handle the case
where $(M, \omega)$ is {\em Futaki nondegenerate}. We end this
section by defining right inverses for the corresponding operator on
$(N_j , \eta_j)$.

\subsection{Analysis of the operators defined on $(M^{*}, \omega)$}

Assume that $(M, \omega)$ is a compact \K\ manifold or \K\ orbifold
with isolated singularities and further assume that $\omega$ has
constant scalar curvature. We first construct a right inverse for
the operator ${\mathbb L}_g$ when $m\geq 3$ and when $(M,\omega)$ is
{\em nondegenerate}. Next, we proceed with the proof of the
corresponding result when $m=2$ and when the kernel of ${\mathbb
L}_g$ is nontrivial but $(M,\omega)$ is {\em Futaki nondegenerate},
i.e. when the linearized Futaki invariant is nondegenerate.

\medskip

The mapping properties of ${\mathbb L}_g$, when defined between
weighted spaces, depends heavily on the choice of the weight
parameter. Recall that, by definition, $\zeta \in {\mathbb R}$ is an
indicial root of ${\mathbb L}_g$ at $p_j$ if there exists some
nontrivial function $v \in {\mathcal C}^{\infty} (\del B_{j,1})$
such
\begin{equation}
{\mathbb L}_g \, (|z|^{\zeta} \, v ) = {\mathcal O} (|z|^{\zeta -
3}) \label{eq:5.1}
\end{equation}
near $0$ (here we have implicitly used the coordinates defined in
Section 2 to parameterize $M$ close to the point $p_j$).

\medskip

Let $\Delta_0$ denote the Laplacian in ${\mathbb C}^{m}$  with its
standard \K\ form. Using (\ref{eq:2.1}), it is easy to check that,
near each $p_j$, (\ref{eq:5.1}) holds for some function $v$ if and
only if
\[
\Delta^{2}_0 \, (|z|^\zeta \, v ) = {\mathcal O} (|z|^{\zeta -3})
\]
Therefore, the set of indicial roots of ${\mathbb L}_g$ at $p_j$  is
equal to the set of indicial roots at the origin for the operator
$\Delta^{2}_0$ defined on ${\mathbb C}^{m} /\Gamma_j$. This later
turns out to be included in ${\mathbb Z} - \{ 5 - 2m, \ldots, -1 \}$
when $m \geq 3$ and is included in ${\mathbb Z}$ when $m=2$ (observe
that the set of indicial roots depends on the group $\Gamma_j$).
Indeed, let $e$ be an eigenfunction of $\Delta_{S^{2m-1}}$ which is
invariant under the action of $\Gamma_j$ and is associated to the
eigenvalue $\gamma \, (2m-2+\gamma)$, where $\gamma \in {\mathbb
N}$, hence
\[
\Delta_{S^{2m-1}} \, e = - \gamma \, (2m-2+ \gamma ) \, e .
\]
Here we identify $S^{2m-1}$ with the unit sphere in ${\mathbb
C}^{m}$. Then
\[
\Delta^{2}_0 \, (|z|^{\zeta} \, e ) = (\zeta -k) (\zeta -k-2) \,
(\zeta -2+2m+k) \, (\zeta -4 + 2m +k) \, |z|^{\zeta -4} \, e .
\]
Therefore, we find that $k$, $k+2$, $2-2m-k$ and $4-2m -k$ are
indicial roots of $\Delta_0^{2}$ at $0$. Since the eigenfunctions of
the Laplacian on the sphere constitute a Hilbert basis of
$L^2(S^{2m-1})$, we have obtained all the indicial roots of
$\Delta_0^{2}$ at the origin.

\medskip

It is clear that the operator
\[
\begin{array}{rclcllll}
L_\delta' : & {{\mathcal C}}^{4, \alpha}_{\delta} (M^{*})
 & \longrightarrow & {{\mathcal C}}^{0,
\alpha}_{\delta - 4} (M^{*})  \\[3mm]
 & \varphi & \longmapsto & {\mathbb L}_g \, \varphi ,
\end{array}
\]
is well defined. It follows from the general theory in \cite{Mel},
where weighted Sobolev spaces are considered instead of weighted
H\"older spaces, and in \cite{Maz}, where the corresponding analysis
in weighted H\"older spaces is performed (see also \cite{Pac-Riv})
that the operator $L_\delta'$ has closed range and is Fredholm,
provided $\delta$ is not an indicial root of ${\mathbb L}_g$ at the
points $p_1, \ldots, p_n$. Under this condition, some duality
argument (in weighted Sobolev spaces) shows that the operator
$L_\delta'$ is surjective if and only if the operator
$L_{4-2m-\delta}'$ is injective. And, still under this assumption
\begin{equation}
\mbox{dim}  \, \mbox{Ker}  \, L_\delta'  =  \mbox{dim} \,
\mbox{Coker} \, L'_{4-2m-\delta} . \label{eq:5.00}
\end{equation}

Using these, we obtain the~:
\begin{prop} Assume that $m \geq 3$, $\delta
\in (4-2m,0)$ and assume that $(M, \omega)$ is {\em nondegenerate}
so that the kernel of ${\mathbb L}_g$ is spanned by the constant
function. Then, the operator
\[
\begin{array}{rclcllll}
L_\delta : & {{\mathcal C}}^{4, \alpha}_{\delta} (M^{*})  \times
{\mathbb R} & \longrightarrow &  {{\mathcal C}}^{0,
\alpha}_{\delta - 4} (M^{*}) \\[3mm]
& (\varphi , \nu ) & \longmapsto & {\mathbb L}_g \, \varphi +  \nu
\end{array}
\]
is surjective and has a one dimensional kernel spanned by the
constant function. \label{pr:5.1}
\end{prop}
{\bf Proof~:} We claim that, when $\delta \in (4-2m, 0)$, the
operator $L_{\delta}'$ has a one dimensional kernel spanned by the
constant function. Indeed, when $\delta \in (4-2m, 0)$, standard
regularity theory implies that the isolated singularities of any
element of the kernel of $L_{\delta}'$  are removable and hence, the
elements of the kernel of $L_{\delta}'$ are in fact a smooth
functions in $M$. Therefore, it follows from our assumption that the
kernel of $L_\delta'$ reduces to the constant functions. It follows
from (\ref{eq:5.00}) that the operator $L_\delta'$ also has a one
dimensional cokernel, which is easily seen to be spanned by the
constant function since (thanks to (\ref{eq:4.3}))
\[
\int_M {\mathbb L}_g  \, \varphi \, dv_g  =0,
\]
for any $\varphi \in {\mathcal C}^{4,\alpha}_\delta (M^{*})$. This
completes the proof of the result. \hfill $\Box$

\medskip

When $m=2$, the above result has to be modified (since $4-2m=0$ in
this case !). We set
\[
{\mathcal D} : = \mbox{Span} \{ \chi_1  , \ldots,  \chi_n \},
\]
where $\chi_j$ is a cutoff function which is identically equal to
$1$ in $B_{j , r_0/2}$ and identically equal to $0$ in $M-B_{j
,r_0}$.  This time, we have the~:
\begin{prop}
Assume that $m=2$, $\delta \in (0,1)$ and assume that $(M, \omega)$
is {\em nondegenerate} so that the kernel of ${\mathbb L}_g$ is
spanned by the constant function. Then
\[
\begin{array}{rclclll}
L_\delta   : &  \left( {{\mathcal C}}^{4, \alpha}_{\delta} (M^{*})
\oplus {\mathcal D}  \right) \times {\mathbb R}& \longrightarrow &
{{\mathcal C}}^{0, \alpha}_{\delta - 4} (M^{*})  \\[3mm]
& (\varphi , \nu ) & \longmapsto & {\mathbb L}_g \, \varphi + \nu
\end{array}
\]
is surjective and has a one dimensional kernel spanned by the
constant function. \label{pr:5.2}
\end{prop}
{\bf Proof~:} We keep the notations of the previous proof. Assume
that $\delta >0$. Then the operator $L_{\delta}'$ is injective
(since we have assumed that the kernel of ${\mathbb L}_g$ is spanned
by the constant function and a nonzero function does not belong to
${{\mathcal C}}^{4, \alpha}_{\delta} (M^{*})$ when $\delta >0$).
Therefore, when $\delta >0$, $\delta \notin {\mathbb N}$, the
operator $L_{-\delta}'$ is surjective and admits a right inverse,
which, unfortunately, is not unique.

\medskip

Moreover, when $\delta \in (0,1)$, a relative index argument
\cite{Mel} shows that the dimension of the kernel of $L_{-\delta}'$
and the dimension of the cokernel of $L_{\delta}'$  are both equal
to $n$. The kernel of $L_{-\delta}'$ is rather explicit since it is
spanned by the constant function and, for $j=1, \ldots, n-1$, the
unique function $\gamma_j$ solution (in the sense of distributions)
of
\[
{\mathbb L}_g \, \gamma_j = \delta_{p_{j+1}} - \delta_{p_j}
\]
and whose mean value over $M$ is $0$.

\medskip

Let us now assume that $\delta \in (0,1)$. Given $\psi \in {\mathcal
C}^{0, \alpha}_\delta (M^{*})$, we choose $\nu \in {\mathbb R}$ to
be equal to the mean value of the function $\psi$. Since
$L_{-\delta}'$ is surjective, we have the existence of a solution of
\[
{\mathbb L}_g \, \varphi  = \psi - \nu ,
\]
which belongs to ${\mathcal C}^{4, \alpha}_{-\delta} (M^{*})$ (this
solution is for example obtained by applying to $\psi-\nu$ a given
right inverse for $L_{-\delta}'$). It follows from elliptic
regularity theory that, near any $p_j$, the function $\varphi$ can
be expanded as
\[
\varphi (z) = a_j + b_j \, \log |z| + \tilde \varphi_j (z)
\]
where $a_j, b_j \in {\mathbb R}$ and $\tilde \varphi_j \in {\mathcal
C}^{4, \alpha}_{\delta} (B_{j, r_0}^{*})$. This implies that the
function $\varphi$ is a solution (in the sense of distributions) of
\begin{equation}
{\mathbb L}_g \, \varphi + \nu = \psi - c_2 \, \sum_{j} b_j \,
\delta_{p_j} \label{eq:5.0}
\end{equation}
where $c_2 = 2 \, |S^{3}| \neq 0$. Using the fact that the functions
$\gamma_j$ are in the kernel of $L_{-\delta}'$, we can assume
without loss of generality that the $b_j$ at the different points
$p_j$ are all equal, by adding to $\varphi$ a suitable linear
combination of the functions $\gamma_j$ (this amounts to choose a
particular right inverse of $L_{-\delta}'$). Integration of
(\ref{eq:5.0}) over $M$ implies that $0 =  - c_2 \, \sum_j b_j$.
Hence all $b_j$ are equal to $0$ and, near $p_j$, the function
$\varphi$ can be expanded as
\[
\varphi (z) = a_j + \tilde \varphi_j (z).
\]
This shows that there exists a choice of the right inverse
$G_{-\delta}'$ of $L_{-\delta}'$ such that, if $\psi \in {\mathcal
C}^{0,\alpha}_{\delta-4}(M^{*})$ and if $\nu$ is the mean value of
$\psi$, then
\[
G_{-\delta}' (\psi-\nu) \in {\mathcal C}^{4, \alpha}_{\delta}(M^{*})
\oplus {\mathcal D}.
\]
This completes the proof of the result.  \hfill $\Box$

\begin{remar}Observe that, given $\psi  \in {\mathcal C}^{0,
\alpha}_{\delta-4} (M^{*})$, the constant $ \nu \in {\mathbb R}$ in
the equation ${\mathbb L}_g \, \varphi + \nu = \psi$ is equal to the
mean value of $\psi$ so that $\psi - \nu $ is $L^{2}$-orthogonal to
the kernel of ${\mathbb L}_g$ which is spanned by the constant
function. \label{re:5.1}
\end{remar}

We turn to the case where the kernel of ${\mathbb L}_g$ is not only
spanned by the constant function and we now assume that $(M,
\omega)$ is {\em Futaki nondegenerate}. The proof relies on the
following result which replaces Proposition~\ref{pr:5.1} and whose
proof is identical.
\begin{prop}
Assume that $m \geq 3$ and that $(M, \omega)$ is {\em Futaki
nondegenerate}. Then, for all $\delta \in (4-2m,0)$ the operator
\[
\begin{array}{rclclll}
L_\delta : & {\mathcal C}^{4, \alpha}_\delta (M^{*}) \times
 {\mathcal H}^{1,1}(M, {\mathbb C}) \times {\mathbb R} & \longrightarrow &
{\mathcal C}^{0,
\alpha}_{\delta -4} (M^{*})\\[3mm]
& (\varphi , \beta,  \nu ) & \longmapsto  & {\mathbb L}_g \, \varphi
+ {\mathcal L}_g \, \beta + \nu ,
\end{array}
\]
is surjective and has a kernel which is equal to the kernel of
${\mathbb L}_g$. \label{pr:5.3}
\end{prop}

Given $\psi  \in {\mathcal C}^{0, \alpha}_{\delta -4} (M^{*})$, the
$(1,1)$ form $\beta \in {\mathcal H}^{1,1}(M, {\mathbb C})$ and the
constant $\nu \in {\mathbb R}$ in the equation ${\mathbb L}_g \,
\varphi + {\mathcal L}_g \, \beta + \nu  =  \psi$ are chosen in such
a way that $\psi - {\mathcal L}_M \, \beta - \nu $ is
$L^{2}$-orthogonal to the elements of the kernel of ${\mathbb L}_g$.

\medskip

Since the space of holomorphic vector fields is finite dimensional
so is $\mbox{Ker} \, {\mathbb L}_g$ and one can replace in the above
statement the space ${\mathcal H}^{1,1} (M, {\mathbb C})$ by a
finite dimensional subspace $D\subset {\mathcal H}^{1,1} (M,
{\mathbb C})$. We claim that this subspace can in turn be replaced
by $D_{\bar r_0}$ a finite dimensional space of closed $(1,1)$ forms
which are supported in $M_{\bar r_0}$, provided $\bar r_0$ is fixed
small enough. Indeed, near each $p_j$, any element of $\beta \in D
\subset {\mathcal H}^{1,1} (M, {\mathbb C})$ can be decomposed as
\[
\beta = d \, \pi_j
\]
We truncate the potential $\pi_j$ between $2 \bar r_0$ and $\bar
r_0$ and define
\[
\beta_{\bar r_0} : = d ( (1-\chi_{\bar r_0}) \, \pi_j),
\]
where $\chi_{\bar r_0}$ is a cutoff function identically equal to
$0$ in $M_{2 \, \bar r_0}$ and identically equal to $1$ in each
$B_{j ,\bar r_0}$. We set
\[
D_{\bar r_0} : =  \mbox{Span} \, \{ \beta_{\bar r_0} \quad : \quad
\beta \in D \}.
\]
When $m \geq 3$, it is easy to check that, given $\delta \in (4-2m,
0)$, the operator
\[
\begin{array}{rclclll}
L_\delta' : & {\mathcal C}^{4, \alpha}_\delta (M^{*}) \times
 D_{\bar r_0} \times {\mathbb R} & \longrightarrow & {\mathcal C}^{0,
\alpha}_{\delta -4} (M^{*})\\[3mm]
& (\varphi, \beta, \nu ) & \longmapsto  & {\mathbb L}_M \, \varphi +
{\mathcal L}_M \, \beta + \nu
\end{array}
\]
is surjective provided $\bar r_0$ is chosen small enough.

\medskip

In dimension $m=2$, this result has to be modified. As above we find
that, given $\delta \in (0,1)$, the operator
\[
\begin{array}{rclclll}
L_\delta' : & ({\mathcal C}^{4, \alpha}_\delta (M^{*}) \oplus
{\mathcal D}) \times D_{\bar r_0} \times {\mathbb R} &
\longrightarrow &
{\mathcal C}^{0, \alpha}_{\delta -4} (M^{*})\\[3mm]
& (\varphi, \beta, \nu ) & \longmapsto  & {\mathbb L}_g \, \varphi +
{\mathcal L}_g \, \beta +  \nu
\end{array}
\]
is surjective and has a kernel which is equal to the kernel of
${\mathbb L}_g$.

\subsection{Operators defined on $(N_j , \eta_j)$}

Assume that $(N_j, \eta_j)$ is an ALE space with zero scalar
curvature \K\ metric $\eta_j$. Further assume that, at infinity, the
\K\ form $\eta_j$ can  be expanded as
\[
\eta_j =  i \, \del \, \bar \del \, ( \mbox{$\frac{1}{2}$} \,
|u|^{2} + \tilde \varphi_{j}(u)),
\]
where $\tilde \varphi_j$ satisfies
\[
\nabla^{2} \, \tilde \varphi_j \in {\mathcal C}^{2, \alpha}_{2-2m}
(C_{j, R_0}).
\]
We denote by $g_j$ the metric associated to the \K\ form $\eta_j$.
Again the analysis of ${\mathbb L}_{g_j}$ when defined between
weighted spaces follows from the general theory developed in
\cite{Mel} and \cite{Maz} (see also \cite{Pac-Riv}) and the mapping
properties of ${\mathbb L}_{g_j}$ when defined between weighted
spaces depends heavily on the choice of the weight parameter.

\medskip

Recall that $\zeta \in {\mathbb R}$ is an indicial root of ${\mathbb
L}_{g_j}$ at infinity if there exists some nontrivial function $v
\in {\mathcal C}^{\infty} (\del B_{j,1})$ such
\begin{equation}
{\mathbb L}_{g_j} \, (|u|^{\zeta} \, v ) = {\mathcal O} (|u|^{\zeta
-5}) \label{eq:5.2}
\end{equation}
near $\infty$ (we have implicitly used the coordinates defined in
Section 2 to parameterize $N_j$ near its end).

\medskip

Again, it is easy to check that, (\ref{eq:5.2}) holds for some
function $v$ if and only if
\[
\Delta^{2}_0 \, (|u|^{\zeta} \, v ) = {\mathcal O} (|u|^{\zeta - 5})
\]
(here one uses the fact that $g_j = g_{eucl} + {\mathcal O}
(|z|^{2-2m})$ at infinity and hence the coefficients of the Ricci
tensor at infinity are bounded by a constant times $|u|^{-2m}$).
Therefore, the set of indicial roots of ${\mathbb L}_{g_j}$ at
infinity is equal to the set of indicial roots at infinity for the
operator $\Delta^{2}_0$ defined on ${\mathbb C}^{m} /\Gamma_j$.
Again, this set is included in ${\mathbb Z} - \{ 5 - 2m, \ldots, -1
\}$ when $m \geq 3$ and is included in ${\mathbb Z}$ when $m=2$ (the
set of indicial roots depends on the group $\Gamma_j$). The proof of
this fact follows the analysis done in Section 5.1.

\medskip

The operator
\[
\begin{array}{rclcllll}
\tilde L_\delta : & {{\mathcal C}}^{4, \alpha}_{\delta} (N_j)
 & \longrightarrow & {{\mathcal C}}^{0,
\alpha}_{\delta - 4} (N_j)  \\[3mm]
 & \varphi & \longmapsto & {\mathbb L}_{g_j} \, \varphi ,
\end{array}
\]
is well defined (again one uses the fact that $g_j = g_{eucl} +
{\mathcal O} (|z|^{2-2m})$ at infinity). Moreover, according to
\cite{Mel} and \cite{Maz} (see also \cite{Pac-Riv}), this operator
has closed range and is Fredholm, provided $\delta$ is not an
indicial root of ${\mathbb L}_{g_j}$ at infinity. Under this
condition, some duality argument (in weighted Sobolev spaces) shows
that the operator $\tilde L_\delta$ is surjective if and only if the
operator $\tilde L_{4-2m-\delta}$ is injective. And, still under
this assumption
\begin{equation}
\mbox{dim}  \, \mbox{Ker}  \, \tilde L_\delta  =  \mbox{dim} \,
\mbox{Coker} \, \tilde L_{4-2m-\delta} \label{eq:5.000}
\end{equation}

\medskip

The construction of a right inverse for the operator ${\mathbb
L}_{g_j}$ relies on the following result whose proof is essentially
borrowed from \cite{ks}~:
\begin{prop}
\label{pr:5.4} Assume that $(N_j, \eta_j)$ is a constant scalar
curvature ALE \K\ manifold or \K\ orbifold with isolated
singularities. Then, there is no nontrivial solution of ${\mathbb
L}_{g_j} \, \varphi = 0$, which belongs to ${\mathcal C}^{4,
\alpha}_{\delta} (N_j)$, for some $\delta <0$.
\end{prop}
{\bf Proof~:} Assume that ${\mathbb L}_{g_j} \, \varphi =0$ and that
$\varphi \in{\mathcal C}^{4, \alpha}_{\delta} (N_j)$, for some
$\delta <0$. Then, as explained in Section 4,  the vector field
$\del^{\#}_{g_j} \varphi$ is a holomorphic vector field which tends
to $0$ at infinity. Indeed, we have
\[
{\mathbb L}_{g_j} =  2 \, (\bar \del \, \del^\#_{g_j})^{\star} \, (
\bar \del \, \del^\#_{g_j})
\]
and, multiplying ${\mathbb L}_{g_j} \,
\varphi = 0$ by $\varphi$ and integrating by parts, we get \[
\int_{N_j} |\bar \del \, \del^\#_{g_j} \, \varphi |^2 \, dv_{g_j} =0
\]
All integrations are justified because of the decaying behavior of
$\varphi$ at infinity which implies that $\varphi \in {\mathcal
C}^{4, \alpha}_{4-2m} (N_j)$ when $m\geq 3$. Therefore
$\del^\#_{g_j} \, \varphi =0$. Using Hartogs' Theorem, the
restriction of $\del^{\#}_{g_j} \varphi $ to $C_{j , R_0}$ can be
extended to a holomorphic vector field on ${\mathbb C}^{m}$. Since
this vector field decays at infinity, it has to be identically equal
to $0$. This implies that $\del^{\#}_{g_j} \varphi$ is identically
equal to $0$ on $C_{j , R_0}$. However $\varphi$ being a real valued
function, this implies that $\del \varphi = \bar \del \varphi = 0$
in $C_{j ,R_0}$. Hence the function $\varphi$ is constant in $C_{j ,
R_0}$ and decays at infinity. This implies that $\varphi$ is
identically equal to $0$ in $C_{j , R_0}$ and satisfies ${\mathbb
L}_{g_j} \, \varphi =0$ in $N_j$. Now, we use the unique
continuation theorem for solutions of linear elliptic equations to
conclude that $\varphi$ is identically equal to $0$ in $N_j$. \hfill
$\Box$

\medskip

This being understood, we have the~:
\begin{prop}
Assume that $\delta \in (0,1)$. Then
\[
\begin{array}{rclcllll}
\tilde L_\delta : & {{\mathcal C}}^{4, \alpha}_{\delta} (N_j) &
\longrightarrow &  {{\mathcal C}}^{0,\alpha}_{\delta - 4} (N_j)  \\[3mm]
 & \varphi & \longmapsto & {\mathbb L}_{g_j} \, \varphi
\end{array}
\]
is surjective and has a one dimensional kernel spanned by the
constant function. \label{pr:5.5}
\end{prop}
{\bf Proof~:} It follows from Proposition~\ref{pr:5.4} that, when $
\delta' < 0 $ the operator $\tilde L_{\delta'}$ is injective and
this implies that $\tilde L_{\delta}$ is surjective whenever $\delta
> 4-2m$ is not an indicial root of ${\mathbb L}_{g_j}$ at
infinity. \hfill $\Box$

\subsection{Bi-harmonic extensions}

We end up this section by the following simple result whose proof
follows at once from the application of the maximum principle. Here,
as usual, $\Gamma$ is a finite subgroup of $U(n)$ acting freely on
${\mathbb C}^{n}-\{0\}$. We define
\[
\bar B_\Gamma : = \{ z \in {\mathbb C}^{n} \, / \, \Gamma \quad :
\quad |z| \leq 1\},
\]
\[
\bar B_\Gamma^{*} : = \{ z \in {\mathbb C}^{n} \, / \, \Gamma \quad
: \quad |z| \leq 1\},
\]
and
\[
\bar C_\Gamma : = \{ z \in {\mathbb C}^{n} \, / \, \Gamma \quad :
\quad |z| \geq  1\}.
\]
Therefore, when $\Gamma= \Gamma_j$, we have $\bar B_{\Gamma_j} =
\bar B_{j ,1}$, $\bar B_{\Gamma_j}^{*} = \bar B_{j ,1}^{*}$ and
$\bar C_{\Gamma_j} = \bar C_{j ,1}$. Recall that $\Delta_0$ denotes
the Laplacian in ${\mathbb C}^{m}$ with the standard \K\ form. With
these notations in mind, we have
\begin{prop}
Assume that $m \geq 3$. Given $h \in {\mathcal C}^{4, \alpha} (\del
B_\Gamma )$, $k \in {\mathcal C}^{2, \alpha} (\del B_\Gamma)$ there
exist bi-harmonic functions $H^{i}_{h,k} \in {\mathcal C}^{4,
\alpha} (\bar B_\Gamma )$ and $H^{o}_{h,k} \in {\mathcal C}^{4 ,
\alpha}_{4-2m} (\bar C_\Gamma)$ such that
\[
\begin{array}{rlll}
\Delta_0^{2} \,  H^{i}_{h,k} & = & 0 \qquad & \mbox{in} \qquad
B_\Gamma  \\[3mm]
\Delta_0^{2} \, H^{o}_{h,k} & = & 0, \qquad & \mbox{in} \qquad
C_\Gamma ,
\end{array}
\]
with
\[
H^{i}_{h,k} = H^{o}_{h,k} = h \quad \mbox{and} \quad \Delta_0
H^{i}_{h,k} = \Delta_0 H^{o}_{h,k} =k \qquad \mbox{on} \qquad \del
B_\Gamma  .
\]
Moreover,
\[
\| H^{i}_{h,k} \|_{{\mathcal C}^{4, \alpha} (\bar B_\Gamma ) } + \|
H^{o}_{h,k}\|_{{\mathcal C}^{4, \alpha} _{4-2m} (\bar C_\Gamma ) }
\leq c \, (\|h\|_{{\mathcal C}^{4, \alpha}(\del B_\Gamma )} + \|k
\|_{{\mathcal C}^{2, \alpha}(\del B_\Gamma )}).
\]
\label{pr:5.6}
\end{prop}

For later use, it will be convenient to get explicit formulas for
$H^{i}_{h,k}$ and $H^{o}_{h,k}$. We decompose both functions $h$ and
$k$ over eigenfunctions of the Laplacian on the sphere. Namely
\[
h = \sum_{\gamma = 0}^\infty h^{(\gamma)} \, e_\gamma \qquad
\mbox{and} \qquad k = \sum_{\gamma=0}^{\infty} k^{(\gamma)} \,
e_\gamma
\]
where the function $e_\gamma$ satisfies
\[
\Delta_{S^{2m-1}} \, e_\gamma = - \gamma \,  (2m -2 + \gamma) \,
e_\gamma
\]
and is normalized so that $\| e_\gamma \|_{L^{2}} =1$. Observe that
we only have to consider the eigenvalues corresponding to
eigenfunctions which are invariant under the action of $\Gamma$.
Then, the functions $H^{i}_{h,k}$ and $H^{o}_{h,k}$ are explicitly
given by
\begin{equation}
H^{i}_{h,k} (z) = \sum_{\gamma=0}^{\infty} \left( \left(
h^{(\gamma)} - \frac{k^{(\gamma)}}{4(m+\gamma)} \right) \,
|z|^{\gamma} + \frac{k^{(\gamma)}}{4(m+\gamma)} \, |z|^{\gamma+2} \,
\right) \, e_\gamma \label{eq:5.44} \end{equation} and
\begin{equation}
H^{o}_{h,k}  (z)= \sum_{\gamma=0}^{\infty} \left( \left(
h^{(\gamma)} + \frac{k^{(\gamma)}}{4(\gamma+m-2)} \, \right) \,
|z|^{2-2m -\gamma} - \frac{k^{(\gamma)}}{4(\gamma+m-2)} \,
|z|^{4-2m-\gamma} \, \right) \, e_\gamma \label{eq:5.444}
\end{equation}

{\bf Proof of Proposition~\ref{pr:5.6}~:} The existence of
$H^{i}_{h,k}$ is clear and the estimate follows at once.

\medskip

The explicit expression of $H^o_{h,k}$ provides a direct proof of
the estimate of this function. First observe that elliptic
regularity implies that, there exists $c = c(m) >0$ and $N = N(m)
\in {\mathbb N}$ such that
\[
\|e_\gamma\|_{L^\infty} \leq c \, (1+\gamma)^{N} \, \|
e_\gamma\|^2_{L^2} = c \, (1+ \gamma)^N
\]
since we have normalized the functions $e_\gamma$ to have $L^2$-norm
equal to $1$. In addition, Cauchy-Schwartz inequality yields
\[ |h^{(\gamma)}|+ |k^{(\gamma)}|\leq c \,( \| h\|_{{\mathcal C}^{4,
\alpha}} +  \| k \|_{{\mathcal C}^{2, \alpha}} )\] for some constant
which does not depend on $\gamma$. Using these two information
together with (\ref{eq:5.444}) we conclude that
\[
\sup_{|z|\geq 2} \left( |z|^{2m-4} |H^o_{h,k}| + |z|^{2m-2} \,
|\Delta_0 \, H^o_{h,k}| \right) \leq c \, ( \| h\|_{{\mathcal C}^{4,
\alpha}} +  \| k \|_{{\mathcal C}^{2, \alpha}} )
\]
since the series are absolutely convergent for $|z|$ larger than
$1$. The maximum principle applied in $\{ z \in C_{\Gamma} \quad :
\quad |z | \in [1,2]\}$ then allows to fill in the gap in the
estimate and we conclude that
\[
\sup_{|z|\geq 1} \left( |z|^{2m-4} |H^o_{h,k}| + |z|^{2m-2} \,
|\Delta_0 \, H^o_{h,k}| \right) \leq c \, ( \| h\|_{{\mathcal C}^{4,
\alpha}} +  \| k \|_{{\mathcal C}^{2, \alpha}} )
\]
The estimates for the derivatives of $H^{o}_{h,k}$ follow from
Schauder's estimates. \hfill $\Box$

\medskip

When $m=2$, the result has to be slightly modified since in this
case we can choose
\begin{equation}
H^{o}_{h,k} (z) = h^{(0)} \, |z|^{-2} + \frac{k^{(0)}}{2} \, \log
|z| + \sum_{\gamma = 1}^{\infty} \left( \left( h^{(\gamma)} +
\frac{k^{(\gamma)}}{4\gamma} \, \right) \, |z|^{-2 - \gamma}  -
 \frac{k^{(\gamma)}}{4\gamma} \,
|z|^{-\gamma} \,  \right) \, e_\gamma \label{eq:5.4444}
\end{equation} This time, one can check that
\begin{equation} H^{o}_{h,k} \in {\mathcal C}^{4,
\alpha}_{-1} (\bar C_\Gamma) \oplus \mbox{Span} \, \{\log |z|\}
\label{eq:5.4}
\end{equation}
and that
\begin{equation}
\| H^{o}_{h,k}\|_{{\mathcal C}^{4, \alpha} _{-1} (\bar C_\Gamma )
 \oplus  \mbox{Span} \, \{\log |z|\}} \leq c \, (\|h\|_{{\mathcal C}^{4, \alpha}(\del B_\Gamma )} + \|k
\|_{{\mathcal C}^{2, \alpha}(\del B_\Gamma )}). \label{eq:5.04}
\end{equation}

\section{Constant scalar curvature \K\ metrics}

We set
\[
r_\e : =  \e^{\frac{n-1}{n}} \qquad \mbox{and} \qquad R_\e : =
\frac{r_\e}{\e}  = \e^{-\frac{1}{n}}.
\]

\subsection{Perturbation of $\omega$}

We will now use the result of the previous sections to perturb
$\omega$, the \K\ form on $M_{r_\e}$, into infinite families of
constant scalar curvature \K\ forms which are defined on $M_{r_\e}$
and which are parameterized by the boundary data of their
potentials. We carry this analysis when $(M, \omega)$ is {\em Futaki
nondegenerate} since this is the most complete case. We consider the
perturbed \K\ form
\begin{equation}
\tilde \omega = \omega + i \, \del \, \bar \del \, \varphi + \beta .
\label{eq:6.1}
\end{equation}
where $\beta$ is a closed $(1,1)$ form and $\varphi$ is a function
defined on $M_{r_\e}$. The scalar curvature of $\tilde \omega$ is
given by
\begin{equation}
{\bf s} (\tilde \omega) = s (\omega) -  ({\mathbb L}_g \, \varphi +
{\mathcal L}_g \, \beta ) + Q_g (\nabla^{2} \varphi , \beta) ,
\label{eq:6.2}
\end{equation}
where the operators ${\mathbb L}_g$ and ${\mathcal L}_g$ have been
defined in (\ref{eq:4.1}) and (\ref{eq:4.2}) and where $Q_g$
collects all the nonlinear terms. The structure of $Q_g$ is quite
complicated however, away from the support of the elements of
$D_{\bar r_0}$ (i.e. in each $\bar B_{j,\bar r_0}$), we have $Q_g
(\nabla^{2} \varphi, \beta) = Q_g (\nabla^{2} \varphi, 0)$ and this
operator, only acting on the function $\varphi$, enjoys the
following decomposition
\begin{equation}
\begin{array}{rllllll}
Q_g ( \nabla^{2} \varphi ,0 ) & = &  \sum_{q} B_{q,4,2}(\nabla^{4}
\varphi, \nabla^{2} \varphi) \, C_{q,4,2} ( \nabla^{2} \varphi) \\[3mm]
& + & \sum_{q} B_{q,3,3}(\nabla^{3} \varphi, \nabla^{3} \varphi) \,
C_{q,3,3} ( \nabla^{2} \varphi) \\[3mm]
& + & |z| \, \sum_{q} B_{q,3,2}(\nabla^{3} \varphi, \nabla^{2}
\varphi) \, C_{q,3,2} ( \nabla^{2} \varphi) \\[3mm]
& + & \sum_{q} B_{q,2,2}(\nabla^{2} \varphi, \nabla^{2} \varphi) \,
C_{q,2,2} ( \nabla^{2} \varphi)
\end{array}
\label{eq:6.3}
\end{equation}
where the sum over $q$ is finite, the operators $(U,V)
\longrightarrow B_{q,a,b} (U,V)$ are bilinear in the entries and
have coefficients which are smooth functions on $\bar B_{j,\bar
r_0}$. The nonlinear operators $W \longrightarrow C_{q,a,b} (W)$
have Taylor expansions (with respect to $W$) whose coefficients are
smooth functions on $\bar B_{j,\bar r_0}$. These facts follow at
once from the expression of the scalar curvature of $\tilde \omega$
in local coordinates as given in (\ref{eq:4.0}).

\medskip

We would like to solve the equation
\begin{equation}
{\bf s} (\tilde \omega) = {\bf s} (\omega) + \nu \label{eq:6.4}
\end{equation} in $M_{r_\e}$, where $\nu \in {\mathbb R}$.

\medskip

We fix a constant $\kappa >0$ (large enough). Assume that we are
given boundary data $h_j \in {\mathcal C}^{4, \alpha}(\del
B_{\Gamma_j})$ and $k_j \in {\mathcal C}^{2, \alpha}(\del
B_{\Gamma_j})$, for $j = 1, \ldots, n$, satisfying
\begin{equation}
\| h_j \|_{{\mathcal C}^{4, \alpha} (\del B_{\Gamma_j})}  \leq
\kappa \, r_\e^{4} \qquad \mbox{and} \qquad  \| k_j\|_{{\mathcal
C}^{2, \alpha} (\del B_{\Gamma_j})} \leq \kappa \, r_\e^{4} .
\label{eq:6.5}
\end{equation}
When $m\geq 3$, we define
\begin{equation}
H_{{\bf h},{\bf k}} : =  \sum_j \chi_j \, H^{o}_{h_j, k_j} (\cdot /
r_\e ) , \label{eq:6.6}
\end{equation}
where we have set
\[
{\bf h} :=(h_1, \ldots, h_n) \qquad \mbox{and} \qquad  {\bf k} : =
(k_1, \ldots, k_n),
\]
and where we recall that the cutoff functions $\chi_j$ are
identically equal to $1$ in $B_{j  , r_0/2}$ and identically equal
to $0$ in $M - B_{j ,r_0}$.

\medskip
When $m=2$, some modifications are necessary. We decompose each
$k_j$ as
\[
k_j = k_j^{(0)} + k_j^{\perp} ,
\]
where $k_j^{(0)}$ is a constant function and $k_j^{\perp}$ has mean
$0$ on $\del B_{\Gamma_j}$. With this decomposition in mind, we
define
\begin{equation}
H_{{\bf h}, {\bf k} }  : =  \displaystyle \sum_j \chi_j \, \left(
H^{o}_{h_j, k_j^{\perp}} (\cdot / r_\e ) + \frac{k_j^{(0)}}{2} \,
\log | \cdot | \right)  .
\label{eq:6.07}
\end{equation}

We replace in (\ref{eq:6.1}) the function $\varphi$ by $H_{{\bf
h},{\bf k}} + \varphi$. Then, (\ref{eq:6.4}) leads to the equation
\begin{equation}
{\mathbb L}_g \, ( H_{{\bf h},{\bf k}} + \varphi) + {\mathcal L}_g
\, \beta + \nu =  Q_g ( H_{{\bf h},{\bf k}} + \varphi , \beta) ,
\label{eq:6.7}
\end{equation}
which we would like to solve in $M_{r_\e}$.

\medskip

\begin{defin}
Given $ \bar r \in (0, r_0/2)$, $k \in {\mathbb N}$, $\alpha \in
(0,1)$ and $\delta \in {\mathbb R}$, the weighted space ${\mathcal
C}^{k, \alpha}_{\delta} (M_{\bar r })$ is defined to be the space of
functions $\varphi \in {{\mathcal C}}^{k, \alpha} (M_{\bar r})$
endowed with the norm
\[
\| \varphi \|_{{{\mathcal C}}^{k, \alpha}_{\delta} (M_{\bar r})}  :
= \| \varphi \|_{{{\mathcal C}}^{k, \alpha} (M_{r_0/2})} + \sum_j \,
\sup_{2\bar r \leq r \leq r_0}  r^{-\delta} \, \| \varphi|_{(B_{j,
r_{0}} -B_{j,\bar r})} (r \, \cdot) \|_{{{\mathcal C} }^{k, \alpha}
( \bar B_{j , 1} - B_{j ,1/2})}
\]\label{de:6.1}
\end{defin}

For each $\bar r\in (0, r_0/2)$, will be convenient to define an
"extension" (linear) operator
\[
{\mathcal E}_{\bar r} : {\mathcal C}^{0, \alpha}_{\delta'} (M_{\bar
r}) \longrightarrow {\mathcal C}^{0, \alpha}_{\delta'} (M^{*}) ,
\]
as follows~:
\begin{itemize}

\item[(i)]  In $M_{\bar r}$, ${\mathcal E}_{\bar r} \, (\psi ) = \psi$, \\

\item[(ii)] in each $B_{j ,\bar r } - B_{j , \bar r/2}$
\[
{\mathcal E}_{\bar r} \, (\psi) (z)  = \displaystyle \frac{ 2 \,
|z| -\bar r  }{\bar r }\,  \psi \left( \bar r \,  \frac{z}{|z|} \right), \\
\]

\item[(iii)] in each $B_{j , \bar  r/2}$,  $ {\mathcal E}_{\bar r } \, (
\psi ) =  0 $.
\end{itemize}

\medskip

It is easy to check that there exists a constant $c = c (\delta')
>0$, independent of ${\bar r} \in (0, r_0/2)$, such that
\begin{equation}
\| {\mathcal E}_{\bar r} (\psi ) \|_{{\mathcal C}^{0,
\alpha}_{\delta'} (M^{*})} \leq \, c \, \| \psi \|_{{\mathcal C}^{0,
\alpha}_{\delta'} (M_{\bar r})} . \label{eq:6.8}
\end{equation}

We fix
\[
\delta \in (4-2m, 5-2m).
\]
With the above notations and definitions, we rephrase the equation
we would like to solve as
\begin{equation}
L_\delta \, (\varphi , \beta, \nu) = {\mathcal E}_{r_\e} \left( Q_g
( H_{{\bf h},{\bf k}} + \varphi  , \beta) - {\mathbb L}_g \, H_{{\bf
h},{\bf k}} \right)  \label{eq:6.9}
\end{equation}
where $\varphi \in {\mathcal C}^{4, \alpha}_{\delta} (M^{*})$, when
$m\geq 3$ and $\varphi \in {\mathcal C}^{4, \alpha}_{\delta} (M^{*})
\oplus {\mathcal D}$ when $m=2$, $\beta \in D_{\bar r_0}$ and $\nu
\in {\mathbb R}$ have to be determined. Observe that any solution of
(\ref{eq:6.9}) is a solution of (\ref{eq:6.7}). The advantage of the
latter versus the former is that we can now make use of the analysis
of Section 6.1 which allows us to find $G_\delta$ a right inverse
for the operator $L_\delta$ and rephrase the solvability of
(\ref{eq:6.9}) as a fixed point problem
\[
(\varphi , \beta, \nu) = {\mathcal N} (\e, {\bf h},{\bf k} ; \varphi
, \beta )
\]
where the nonlinear operator ${\mathcal N}$ is defined by
\[
{\mathcal N} (\e , {\bf h},{\bf k} ; \varphi  , \beta ) : = G_\delta
\left( {\mathcal E}_{r_\e} \left( Q_g ( H_{{\bf h},{\bf k}} +
\varphi  , \beta) - {\mathbb L}_g \, H_{{\bf h},{\bf k}} \right)
\right).
\]

It will be convenient to denote
\[
{\mathcal F} : = {\mathcal C}^{4, \alpha}_\delta (M^{*}) \times
D_{\bar r_0} \times {\R}
\]
when $m\geq 3$ and
\[
{\mathcal F} : = ({\mathcal C}^{4, \alpha}_\delta (M^{*}) \oplus
{\mathcal D} )\times D_{\bar r_0} \times {\R}
\]
when $m=2$. This space is naturally endowed with the product norm.

\medskip

We first estimate the terms on the right hand side of (\ref{eq:6.9})
when $\varphi=0$ and $\beta =0$ and next show that ${\mathcal N}$ is
a contraction from a suitable small ball in ${\mathcal F}$. This is
the content of the~:
\begin{lemma}
There exists $c_\kappa = c (\kappa) >0$, $\tilde c_\kappa =  \tilde
c(\kappa) >0$ and there exists $\e_\kappa = \e(\kappa) >0$ such
that, for all $\e \in (0, \e_\kappa)$
\begin{equation}
\| {\mathcal N} ( \e, {\bf h},{\bf k} ; 0 , 0) \|_{{\mathcal F}}
\leq c_\kappa \, r_\e^{2m} . \label{eq:6.10} \end{equation} In
addition,
\begin{equation} \| {\mathcal N} (\e, {\bf h},{\bf k} ;
\varphi, \beta) - {\mathcal N} (\e, {\bf h}, {\bf k} ;  \varphi' ,
\beta' ) \|_{{\mathcal F}} \leq \tilde c_\kappa \, r_\e^{2} \, \|
(\varphi - \varphi', \beta - \beta') \|_{{\mathcal F}}
\label{eq:6.11} \end{equation} and
\begin{equation}
\| {\mathcal N} (\e, {\bf h},{\bf k} ; \varphi, \beta) - {\mathcal
N} (\e, {\bf h}', {\bf k}' ;  \varphi , \beta ) \|_{{\mathcal F}}
\leq \tilde c_\kappa \, r_\e^{2m-4} \, \| ({\bf h}-{\bf h'}, {\bf k}
- {\bf k'})\|_{( {\mathcal C}^{4, \alpha})^n \times ({\mathcal
C}^{2, \alpha})^n} \label{eq:6.12} \end{equation} provided
$(\varphi, \beta, 0) , ( \varphi' , \beta' , 0) \in {\mathcal F}$
satisfy
\[
\| (\varphi, \beta, 0)\|_{{\mathcal F}}   \leq \, 2 \, c_\kappa \,
r_\e^{2m} \qquad \mbox{and} \qquad  \| (\varphi' , \beta', 0 )
\|_{{\mathcal F}} \leq 2 \, c_\kappa \, r_\e^{2m} ,
\]
and ${\bf h} : = (h_1, \ldots, h_n), {\bf h}': =  (h'_1, \ldots,
h'_n), {\bf k} : = (k_1, \ldots, k_n) $ and  ${\bf k}' : = (k'_1,
\ldots , k'_n)$ satisfy (\ref{eq:6.5}). \label{le:6.1}
\end{lemma}
{\bf Proof~:} We give a precise proof of the first estimate. The
other estimates follow from similar considerations. In the proof,
the constants $c^{(\ell)}_\kappa >0$ only depend on $\kappa$.

\medskip

{\bf Step 1} We first estimate ${\mathbb L}_g \, H_{{\bf h},{\bf
k}}$. Using the result of Proposition~\ref{pr:5.6}, together with
(\ref{eq:6.5}), we obtain
\begin{equation}
\| \nabla^{2} H_{{\bf h},{\bf k}} \|_{{\mathcal C}^{2,
\alpha}_{2-2m} (M_{r_\e})} \leq c_\kappa^{(1)} \, r_\e^{2m} .
\label{eq:6.13}
\end{equation}

\medskip

Now observe that, by construction,  $\nabla^{2} H_{{\bf h},{\bf k}}
= 0$ in $M_{r_0}$ and hence ${\mathbb L}_g \,H_{{\bf h},{\bf k}} =0$
in this set. Next, \[ \Delta_0^{2} \, H_{{\bf h},{\bf k}} =0
\] in
each $B_{j ,r_0/2} - B_{j  , r_\e}$, hence \[ {\mathbb L}_g \,
H_{{\bf h},{\bf k}} = \left( {\mathbb L}_g - \frac{1}{2}
\Delta_0^{2} \right) \, H_{{\bf h},{\bf k}}
\]
in this set. Using (\ref{eq:2.1}), we conclude that
\[
\| {\mathbb L}_g \, H_{{\bf h},{\bf k}}  \|_{{\mathcal C}^{0,
\alpha}_{\delta - 4}(M_{r_\e})} \leq c_\kappa^{(2)} \, r_\e^{2m}
\]
and
\[
\int_{M} | {\mathcal E}_{r_\e} \,({\mathbb L}_g  \, H_{{\bf h},{\bf
k}}) |  \, dv_g \leq c_\kappa^{(2)} \, r_\e^{2m}
\]
These two estimates together with the properties of $G_\delta$
immediately imply that
\[
\| G_\delta \, ({\mathcal E}_{r_\e} \,({\mathbb L}_g  \, H_{{\bf
h},{\bf k}})) \|_{{\mathcal F}} \leq c_\kappa^{(3)} \, r_\e^{2m}.
\]

{\bf Step 2} We turn to the estimate of $Q_g (\nabla^{2} H_{{\bf
h},{\bf k}} , 0)$. To this aim, we use the structure of $Q_g$ as
described in (\ref{eq:6.3}) together with (\ref{eq:6.13}) to get
\[
\|Q_g(\nabla^{2} H_{{\bf h},{\bf k}}, 0)\|_{{\mathcal C}^{0,
\alpha}(M_{\bar r_0/2})} \leq c_\kappa^{(4)} \, r_\e^{4m} ,
\]
and
\[
\| {\mathcal E}_{r_\e} \, (B_{q,a,b} (\nabla^{2+a} H_{{\bf h},{\bf
k}}, \nabla^{2+b} H_{{\bf h},{\bf k}}) \, C_{q,a,b} (\nabla^{2}
H_{{\bf h},{\bf k}}) ) \|_{{\mathcal C}^{0, \alpha}_{\delta-4} (\bar
B_{j, \bar r_0})} \leq c_\kappa^{(4)} \, r_\e^{8-a-b -\delta} ,
\]
Therefore, we conclude that
\[
\|{\mathcal E}_{r_\e}  \, ( Q_g (\nabla^{2} H_{{\bf h},{\bf k}} ,0)
)\|_{{\mathcal C}^{0, \alpha}_{\delta -4}(M^{*})} \leq
c_\kappa^{(5)} \, r_\e^{6 -\delta }
\]
as well as
\[
\int_M |{\mathcal E}_{r_\e}  \, ( Q_g (\nabla^{2} H_{{\bf h},{\bf
k}} ,0) )| \, dv_g \leq c_\kappa^{(5)} \, r_\e^{2m+2}
\]
The properties of $G_\delta$ yield
\[
\| G_\delta \, ( {\mathcal E}_{r_\e}  \, ( Q_g (\nabla^{2} H_{{\bf
h},{\bf k}} ,0) ) ) \|_{\mathcal F} \leq c_\kappa^{(6)} \,
r_\e^{6-\delta}
\]
This completes the proof of the first estimate.

\medskip

{\bf Step 3} We now turn to the derivation of the second estimate.
Again, we use the structure of $Q_g$ as described in (\ref{eq:6.3})
to get
\[
\|Q_g(\nabla^{2} H_{{\bf h},{\bf k}} + \varphi , \beta) -
Q_g(\nabla^{2} H_{{\bf h},{\bf k}} + \varphi' , \beta')\|_{{\mathcal
C}^{0, \alpha}(M_{\bar r_0})} \leq c_\kappa^{(7)} \, r_\e^{2m} \, \|
(\varphi- \varphi', \beta - \beta', 0)\|_{\mathcal F} ,
\]
and, arguing as above, we find that
\[
\|{\mathcal E}_{r_\e}  \, ( Q_g (\nabla^{2} H_{{\bf h},{\bf k}} +
\varphi ,\beta ) -Q_g (\nabla^{2} H_{{\bf h},{\bf k}} + \varphi'
,\beta' ) )\|_{{\mathcal C}^{0, \alpha}_{\delta -4}(\bar B_{j, \bar
r_0})} \leq c_\kappa^{(7)} \, r_\e^{2} \, \| (\varphi- \varphi', 0,
0)\|_{\mathcal F}
\]
Therefore, we conclude that
\[
\|{\mathcal E}_{r_\e}  \, ( Q_g (\nabla^{2} H_{{\bf h},{\bf k}} +
\varphi ,\beta ) -Q_g (\nabla^{2} H_{{\bf h},{\bf k}} + \varphi'
,\beta' ) )\|_{{\mathcal C}^{0, \alpha}_{\delta -4}(M^{*})} \leq
c_\kappa^{(8)} \,  r_\e^{2} \, \| (\varphi- \varphi', \beta -
\beta', 0)\|_{\mathcal F}
\]
as well as
\[
\int_M |{\mathcal E}_{r_\e}  \, ( Q_g (\nabla^{2} H_{{\bf h},{\bf
k}} + \varphi ,\beta ) - Q_g (\nabla^{2} H_{{\bf h},{\bf k}} +
\varphi' ,\beta' ))| \, dv_g \leq c_\kappa^{(8)} \,
r_\e^{2m-2+\delta} \, \| (\varphi- \varphi', \beta - \beta',
0)\|_{\mathcal F}
\]
Observe that, in order to derive the second estimate, we have
implicitly used the fact that the computation of the scalar
curvature only involves second and higher partial differential of
the functions $\varphi$ and $\varphi'$ and hence, in dimension
$m=2$,  the effect of the elements of ${\mathcal D}$ have no
influence in $\bar B_{j, r_0} -B_{j, r_\e}$. The estimate then
follows from the boundedness of $G_\delta$.

\medskip

{\bf Step 4} In order to prove the third estimate, we first observe
that
\[ \| {\mathbb L}_g \, (H_{{\bf h},{\bf k}} -H_{{\bf h}',{\bf
k}'} ) \|_{{\mathcal C}^{0, \alpha}_{\delta - 4}(M_{r_\e})} \leq
c_\kappa^{(9)} \, r_\e^{2m-4} \, \| ({\bf h}-{\bf h'}, {\bf k} -
{\bf k'})\|_{( {\mathcal C}^{4, \alpha})^{n} \times ({\mathcal
C}^{2, \alpha})^{n}}
\]
and
\[
\int_{M} | {\mathcal E}_{r_\e} \,({\mathbb L}_g  \, (H_{{\bf h},{\bf
k}}) -H_{{\bf h}',{\bf k}'})) | \, dv_g \leq c_\kappa^{(9)} \,
r_\e^{2m-4} \, \| ({\bf h}-{\bf h'}, {\bf k} - {\bf k'})\|_{(
{\mathcal C}^{4, \alpha})^{n} \times ({\mathcal C}^{2, \alpha})^{n}}
\]
Next, we have
\[
\begin{array}{rllllll}
\|{\mathcal E}_{r_\e}  \, ( Q_g (\nabla^{2} H_{{\bf h},{\bf k}} +
\varphi ,\beta ) -Q_g (\nabla^{2} H_{{\bf h},{\bf k}} + \varphi
,\beta
) )\|_{{\mathcal C}^{0, \alpha}_{\delta -4}(M^{*})} \qquad \qquad \\[3mm]
\leq c_\kappa^{(10)} \,  r_\e^{2-\delta} \, \| ({\bf h}-{\bf h'},
{\bf k} - {\bf k'})\|_{( {\mathcal C}^{4, \alpha})^{n} \times
({\mathcal C}^{2, \alpha})^{n}} \end{array}
\]
as well as
\[
\begin{array}{rlllll}
\displaystyle \int_M |{\mathcal E}_{r_\e}  \, ( Q_g (\nabla^{2}
H_{{\bf h},{\bf k}} + \varphi ,\beta ) - Q_g (\nabla^{2} H_{{\bf
h}',{\bf k}'}
+ \varphi ,\beta ))| \, dv_g \qquad \qquad \\[3mm]
\leq c_\kappa^{(10)} \, r_\e^{2m-2 } \, \| ({\bf h}-{\bf h'}, {\bf
k} - {\bf k'})\|_{( {\mathcal C}^{4, \alpha})^{n} \times ({\mathcal
C}^{2, \alpha})^{n}} \end{array}
\]
The third estimate now follows from the boundedness of $G_\delta$.

\medskip

This completes the proof of the result. \hfill $\Box$

\medskip

Reducing $\e_\kappa >0$ if necessary, we can assume that,
\begin{equation} \tilde c_\kappa \, r_\e^{2}  \leq \frac{1}{2}
\label{eq:6.1400}
\end{equation} for all $\e \in (0, \e_\kappa )$. Then, the estimates
(\ref{eq:6.10}) and (\ref{eq:6.11}) in the above Lemma are enough to
show that
\[ (\varphi, \beta, \nu) \longmapsto {\mathcal N} (\e, {\bf h} ,
{\bf k} ; \varphi , \beta)
\]
is a contraction from
\[
\{ (\varphi, \beta, \nu) \in {\mathcal F} \qquad :  \qquad \|(
\varphi , \beta, \nu) \|_{{\mathcal F}} \leq 2 \, c_\kappa \,
r_\e^{2m} \}
\]
into itself and hence has a unique fixed point $(\varphi_{\e, {\bf
h}, {\bf k}}, \beta_{\e, {\bf h}, {\bf k}} , \nu_{\e, {\bf h}, {\bf
k}})$ in this set. This fixed point is a solution of (\ref{eq:6.7})
in $M_{r_\e}$ and hence provides a constant scalar curvature \K\
form on $M_{r_\e}$.

\medskip

To summarize, we have obtained the~:
\begin{prop}
Given $\kappa >0$, there exists $\hat c_\kappa >0$ and $\e_\kappa
>0$ such that, for all $\e \in (0, \e_\kappa)$, for all $h_j \in
{\mathcal C}^{4, \alpha}(\del B_{\Gamma_j})$ and all $k_j \in
{\mathcal C}^{2, \alpha}(\del B_{\Gamma_j})$ satisfying
(\ref{eq:6.5}), the \K\ form
\[
\omega_{\e, {\bf h}, {\bf k}} : = \omega + i \, \del \, \bar \del \,
\varphi_{\e, {\bf h}, {\bf k}} + \beta_{\e , {\bf h}, {\bf k}} ,
\]
defined on $M_{r_\e}$, has constant scalar curvature equal to
\[
{\bf s} (\omega_{\e, {\bf h}, {\bf k}} ) =  {\bf s} (\omega) +
\nu_{\e, {\bf h}, {\bf k}}.
\]
Moreover, $\beta_{\e, {\bf h}, {\bf k}} \in D_{\bar r_0}$,
\[
\| \varphi_{\e, {\bf h}, {\bf k}} \, |_{B_{j  ,2 r_\e} -B_{j ,r_\e}}
(r_\e \, \cdot) - H_{h_j, k_j}^{o} \|_{{\mathcal C}^{4, \alpha} (
B_{j, 2} - B_{j , 1} )}\leq \hat c_\kappa \, r_\e^{2m+ \delta} ,
\]
and
\[
|  \nu_{\e, {\bf h}, {\bf k} } |\leq \hat c_\kappa \, r_\e^{2m} .
\]
\label{pr:6.1}
\end{prop}

Using (\ref{eq:6.11}) and (\ref{eq:6.12}), and increasing $\hat
c_\kappa$ if this is necessary, one can check that
\begin{equation}
\begin{array}{rllll} \| ( \varphi_{\e, {\bf h}, {\bf k}}  -
\varphi_{\e , {\bf h'}, {\bf k'}})\, |_{B_{j ,2 r_\e} - B_{j ,r_\e}}
(r_\e \, \cdot) & - & H_{h_j -h'_j, k_j - k'_j}^{o}
\|_{{\mathcal C}^{4, \alpha} ( B_{j , 2} - B_{j , 1} )}  \\[3mm]
& \leq & \hat c_\kappa \, r_\e^{2m -4 + \delta} \, \| ({\bf h}-{\bf
h'}, {\bf k}-{\bf k'})\|_{( {\mathcal C}^{4, \alpha})^{n} \times
({\mathcal C}^{2, \alpha})^{n} },
\end{array}
\label{eq:6.14}
\end{equation}
and
\begin{equation}
| \nu_{\e, {\bf h}, {\bf k}} - \nu_{\e, {\bf h'}, {\bf k'}} |\leq
\hat c_\kappa \, r_\e^{2m-4} \, \| ({\bf h}-{\bf h'}, {\bf k} - {\bf
k'})\|_{( {\mathcal C}^{4, \alpha})^{n} \times ({\mathcal C}^{2,
\alpha})^{n}} . \label{eq:6.15}
\end{equation}
Indeed, if
\[
(\varphi , \beta , \nu )  =  {\mathcal N} (\e, {\bf h}, {\bf k} ;
\varphi , \beta) \qquad \mbox{and} \qquad (\varphi' , \beta' , \nu'
) = {\mathcal N} (\e, {\bf h}', {\bf k}' ; \varphi' , \beta') ,
\]
we can write
\[
\begin{array}{rlcllllll}
(\varphi' - \varphi, \beta'  - \beta , \nu' -\nu ) & =  &( {\mathcal
N} (\e, {\bf h}', {\bf k}' ; \varphi' , \beta') - {\mathcal N} (\e,
{\bf h}', {\bf k}' ; \varphi, \beta) ) \\[3mm]
& + & ({\mathcal N} (\e, {\bf h}', {\bf k}' ; \varphi , \beta) -
{\mathcal N} (\e, {\bf h}, {\bf k} ; \varphi , \beta) )
\end{array}
\]
Using (\ref{eq:6.11}) we get \[ \|(\varphi' - \varphi, \beta' -
\beta , \nu' -\nu )\|_{\mathcal F} \leq 2 \, \| ({\mathcal N} (\e,
{\bf h}', {\bf k}' ; \varphi , \beta) - {\mathcal N} (\e, {\bf h},
{\bf k} ; \varphi , \beta) )\|_{\mathcal F}
\]
and the result follows from (\ref{eq:6.12}).

\subsection{Perturbation of $\eta_{j}$}

Now, we would like to perturb the \K\ form on $N_{j,R_\e}$ into some
infinite dimensional family of constant scalar curvature \K\ forms
which are parameterized by their scalar curvature and the boundary
data of their potentials.

\medskip

We consider the perturbed \K\ form
\begin{equation}
\tilde \eta_j  = \eta_j+  i \, \del \, \bar \del \, \varphi ,
\label{eq:6.155} \end{equation} The scalar curvature of $\tilde
\eta_j$ is given by
\begin{equation}
{\bf s} (\tilde \eta_j) =  - {\mathbb L}_{g_j} \, \varphi + Q_{g_j}
(\nabla^{2} \varphi ) , \label{eq:6.16}
\end{equation}
since the scalar curvature of $\eta_j$ is identically equal to $0$.
Again, the structure of the nonlinear operator $Q_{g_j}$ is quite
complicated but, in $C_{j ,R_0}$, it enjoys a decomposition similar
to the one described in (\ref{eq:6.3}). Indeed, using
(\ref{eq:3.3}), (\ref{eq:3.4})-(\ref{eq:3.4bis}) and (\ref{eq:4.0}),
we see that we can decompose
\[
\begin{array}{rlllll}
Q_{g_j} ( \nabla^{2} \varphi ) & = & \sum_{q} B_{q,4,2}(\nabla^{4}
\varphi, \nabla^{2} \varphi) \, C_{q,4,2} ( \nabla^{2} \varphi) \\[3mm]
& + & \sum_{q} B_{q,3,3}(\nabla^{3} \varphi, \nabla^{3} \varphi) \,
C_{q,3,3} ( \nabla^{2} \varphi) \\[3mm]
& + & \sum_{q} |u|^{1-2m} \, B_{q,3,2}(\nabla^{3} \varphi,
\nabla^{2} \varphi) \, C_{q,3,2} ( \nabla^{2} \varphi)
\\[3mm]
& +& \sum_{q} |u|^{-2m} \, B_{q,2,2}(\nabla^{2} \varphi, \nabla^{2}
\varphi) \, C_{q,2,2} ( \nabla^{2} \varphi)
\end{array}
\]
where the sum over $q$ is finite, the operators $(U,V)
\longrightarrow B_{q,a,b} (U,V)$ are bilinear in the entries and
have coefficients which are bounded functions in ${\mathcal C}^{0,
\alpha} (\bar C_{j,R_0})$. The nonlinear operators $W
\longrightarrow C_{q,a,b} (W)$ have Taylor expansion (with respect
to $W$) whose coefficients are bounded functions on ${\mathcal
C}^{0, \alpha} (\bar C_{j,R_0})$. Even though these operators do
depend on $j$ we have not made this dependence apparent in the
notation.

\medskip

We would like to solve the equation \begin{equation} {\bf s} \,
(\tilde \eta _j ) = \e^{2} \, \nu \label{eq:6.17} \end{equation} in
$N_{j,R_\e}$, where $\nu \in {\mathbb R}$ and where we recall that
\[
R_\e : = \frac{r_\e}{\e}.
\]

We fix a constant $\kappa >0$ large enough and assume that we are
given $\nu \in {\mathbb R}$ and  boundary data $h \in {\mathcal
C}^{4, \alpha}(\del B_{\Gamma_j})$ and $k \in {\mathcal C}^{2,
\alpha}(\del B_{\Gamma_j})$ satisfying
\begin{equation}
|\nu |\leq |{\bf s} (\omega)| +1 \qquad \qquad \| h \|_{{\mathcal
C}^{4, \alpha} (\del B_{\Gamma_j})} \leq \kappa \, R_\e^{4-2m},
\qquad \qquad  \| k \|_{{\mathcal C}^{2, \alpha} (\del
B_{\Gamma_j})} \leq \kappa \, R_\e^{4-2m}, \label{eq:6.18}
\end{equation}
We decompose
\[
h = h^{(0)} + h^{\perp}
\]
where $h^{(0)}$ is a constant function and $h^{\perp}$ has mean $0$
on $\del B_{\Gamma_j}$, and we define
\begin{equation}
\begin{array}{rlllll}
\tilde H_{h , k } & : = &  \tilde \chi_j \, H^{i}_{ h^{\perp} , k}
(\cdot / R_\e ) + h^{(0)} \\[3mm]
& = & \tilde \chi_j \, (H^{i}_{h , k} (\cdot / R_\e) - H^{i}_{h , k}
(0) ) + H^{i}_{h , k} (0),
\end{array}
\label{eq:6.19}
\end{equation}
where $\tilde \chi_j$ is a cutoff function which is identically
equal to $1$ in $C_{j , 2R_0}$ and identically equal to $0$ in $N_{j
,R_0}$.

\medskip

Replacing in (\ref{eq:6.155}) the function $\varphi$ by $\tilde H_{h
, k}+ \varphi$, we see that (\ref{eq:6.15}) can be written as
\begin{equation}
{\mathbb L}_{g_j} \, ( \tilde H_{h, k} + \varphi )  = Q_{g_j} (
\nabla^{2} \tilde H_{h , k} + \varphi ) - \e^{2} \, \nu,
\label{eq:6.99}
\end{equation}
which we would like to solve in $N_{j, R_\e}$. Here $\varphi \in
{\mathcal C}^{4, \alpha}_{\delta} (N_j) $ for some $\delta \in
{\mathbb R}$ has to be determined.

\medskip

\begin{defin}
Given $ \bar R > 2 R_0$, $k \in {\mathbb N}$, $\alpha \in (0,1)$ and
$\delta \in {\mathbb R}$, the weighted space ${\mathcal C}^{k,
\alpha}_{\delta} (N_{j , \bar R})$ is defined to be the space of
functions $\varphi  \in {{\mathcal C}}^{k, \alpha} (N_{j , \bar R})$
endowed with the norm
\[
\| \varphi \|_{{{\mathcal C}}^{k, \alpha}_{\delta} (N_{j,\bar R})} :
= \| \varphi \|_{{{\mathcal C}}^{k, \alpha} (N_{j, 2 R_0})} +
\sup_{2R_0 \leq R \leq \bar R}  R^{-\delta} \, \| \varphi |_{(\bar
C_{j, R_0} -C_{j, \bar R})} (R \, \cdot) \|_{{{\mathcal C} }^{k,
\alpha} ( \bar B_{j , 1} - B_{j ,1/2})}
\]\label{de:6.2}
\end{defin}

For each $\bar R \geq 2 \, R_0$, will be convenient to define an
"extension" (linear) operator
\[
\tilde {\mathcal E}_{\bar R} : {\mathcal C}^{0, \alpha}_{\delta'}
(N_{j , \bar R}) \longrightarrow {\mathcal C}^{0, \alpha}_{\delta'}
(N_j ) ,
\]
as follows~:
\begin{itemize}

\item[(i)] In $N_{j  , R_0}$,  $\tilde {\mathcal E}_{\bar R} \, (\psi ) =
\psi$,\\

\item[(ii)] in $C_{j ,{\bar R}} - C_{j ,2 {\bar R} }$
\[ \tilde {\mathcal E}_{\bar R} \, (\psi ) (u) = \frac{2 \, {\bar R} -
|u|}{\bar R} \, \psi \left(\bar R \, \frac{u}{|u|}\right),
\]

\item[(iii)]  in $C_{j , 2 \, \bar R }$, $\tilde {\mathcal E}_{\bar R} \,
(\psi )=0$.

\end{itemize}
It is easy to check that there exists a constant $c = c( \delta' )
>0$, independent of $\bar R \geq 2 \, R_0$, such that
\begin{equation}
\| \tilde {\mathcal E}_{\bar R} (\psi ) \|_{{\mathcal C}^{0,
\alpha}_{\delta'} (N_j)} \leq \, c \, \| \psi \|_{{\mathcal C}^{0,
\alpha}_{\delta'} (N_{j , {\bar R}})} . \label{eq:6.20}
\end{equation}

We fix
\[
\delta \in (0,1)
\]
The equation we would like to solve can be rewritten as
\begin{equation}
\tilde L_\delta \, \varphi = \tilde {\mathcal E}_{R_\e} \left(
Q_{g_j} ( \tilde H_{h , k}  + \varphi ) - {\mathbb L}_{g_j} \,
\tilde H_{h , k} - \e^{2} \, \nu \right) . \label{eq:6.21}
\end{equation}
where $\varphi \in {\mathcal C}^{4, \alpha}_{\delta} (N_j)$ has to
be determined. Observe that any solution of (\ref{eq:6.21}) is a
solution of (\ref{eq:6.19}). Again, we make use of the analysis of
Section 6.2 in order to find $\tilde G_\delta$ a right inverse for
the operator $\tilde L_\delta$ and rephrase the solvability of
(\ref{eq:6.21}) as a fixed point problem.
\begin{equation}
\varphi = \tilde {\mathcal N}_j (\e, h, k , \nu ; \varphi)
\label{eq:6.22} \end{equation} where the nonlinear operator ${\tilde
{\mathcal N}}$ is defined by
\[
\tilde {\mathcal N} (\e , h, k, \nu  ; \varphi ) : = \tilde G_\delta
\, \left( \tilde {\mathcal E}_{R_\e} \left( Q_{g_j} ( \tilde H_{h,k}
) - {\mathbb L}_{g_j} \, \tilde H_{h, k} - \e^{2} \, \nu \right)
\right)
\]

To keep notations short, it will be convenient to define
\[
\tilde {\mathcal F} : = {\mathcal C}^{4, \alpha}_\delta (N_j)
\]
We first estimate the terms on the right hand side of
(\ref{eq:6.22}) when $\varphi=0$ and next show that $\tilde{\mathcal
N}$ is a contraction from a suitable small ball in $\tilde {\mathcal
F}$. This is the content of the~:
\begin{lemma}
There exists $c >0$ (independent of $\kappa$),  $\tilde c_\kappa =
\tilde c(\kappa) >0$ and there exists $\e_\kappa = \e(\kappa) >0$
such that, for all $\e \in (0, \e_\kappa)$
\begin{equation}
\| \tilde {\mathcal N} ( \e,  h, k, \nu ; 0 ) \|_{\tilde {\mathcal
F}} \leq c \, R_\e^{4-2m -\delta} , \label{eq:6.23}
\end{equation} Moreover, for all $\varphi , \varphi' \in \tilde
{\mathcal F}$, satisfying
\[
\|\varphi \|_{\tilde {\mathcal F}} \leq 2 \, c \, R^{4-2m - \delta}
\qquad \qquad \|\varphi'\|_{\tilde {\mathcal F}} \leq 2 \, c \,
R_\e^{4-2m -\delta} ,
\]
we have
\begin{equation}
\| \tilde {\mathcal N} ( \e, h, k, \nu ; \varphi ) -\tilde {\mathcal
N} ( \e, h, k, \nu ; \varphi' ) \|_{\tilde {\mathcal F}} \leq \tilde
c_\kappa \, R_\e^{4-2m -\delta} \, \| \varphi - \varphi' \|_{\tilde
{\mathcal F}} \label{eq:6.24} \end{equation} and
\begin{equation}
\| \tilde {\mathcal N} ( \e, h, k, \nu ; \varphi ) - \tilde
{\mathcal N} ( \e, h', k', \nu' ; \varphi ) \|_{\tilde {\mathcal F}}
\leq \tilde c_\kappa \, ( R_\e^{-1} \, \| (h-h', k-k')\|_{\mathcal
C^{4, \alpha} \times {\mathcal C}^{2, \alpha}} + R_\e^{4-2m -\delta}
\, |\nu' - \nu|) \label{eq:6.25}
\end{equation} provided $h, h'$ and $k,k'$ satisfy (\ref{eq:6.18}).
\label{le:6.2}
\end{lemma}
{\bf Proof~:} The proof is identical to the proof of
Lemma~\ref{le:6.1}. We give details about the derivation of the
first estimate and leave the two other estimates to the reader

\medskip

It follows from the analysis of Section 5.3, together with
(\ref{eq:6.18}) that
\begin{equation}
\|\nabla^{2} \, \tilde H_{h, k} \|_{{\mathcal C}^{2,
\alpha}_{0}(N_{j , R_\e} )} \leq c_\kappa^{(1)} \, R_\e^{2-2m}
\label{eq:6.26}
\end{equation}
and also that
\begin{equation}
\|\nabla^{2} \, \tilde H_{h, k} \|_{{\mathcal C}^{2,
\alpha}_{0}(\bar C_{j , 2R_0} - C_{j, R_0})} \leq c_\kappa^{(1)} \,
R_\e^{3-2m} \label{eq:6.27}
\end{equation}

We use the fact that, in $C_{j, 2R_0} - C_{j, R_\e}$, we can write
\[
{\mathbb L}_{g_j} \, H_{h, k} = \left( {\mathbb L}_{g_j} -
\frac{1}{2}\, \Delta_0^{2} \right) \, \tilde H_{h , k}.
\]
Then, (\ref{eq:3.4})-(\ref{eq:3.4bis}) together with (\ref{eq:6.26})
yields
\[
\| {\mathbb L}_{g_j} \, \tilde H_{h, k} \|_{{\mathcal C}^{0,
\alpha}_{\delta -4}(N_{j, R_\e})} \leq \, c \, R_\e^{3-2m}
\]

Next, we use the structure of $Q_{g_j}$ together with
(\ref{eq:6.26}) to estimate
\[
\| \tilde {\mathcal E}_{R_\e} \, ( Q_{g_j} ( \tilde H_{h , k} ))
\|_{{\mathcal C}^{0, \alpha}_{\delta -4}(N_j)} \leq c_\kappa^{(2)}
\, R_\e^{6-4m}.
\]
Finally, we estimate
\[
\| \tilde {\mathcal E}_{R_\e} \left( \e^{2} \, \nu \right)
\|_{{\mathcal C}^{0, \alpha}_{\delta-4} (N_j )} \leq \tilde c \,
R_\e^{4-2m-\delta}
\]
for some constant $\tilde c >0$ which does not depend on $\e$ since
$|\nu| \leq  1 +  \, |{\bf s}(\omega )|$.  This completes the proof
of the estimate. \hfill $\Box$

\medskip

Reducing $\e_\kappa >0$ if necessary, we con assume that,
\begin{equation}
\tilde c_\kappa \, R_\e^{4-2m-\delta}  \leq \frac{1}{2}
\label{eq:6.28}
\end{equation} for all $\e \in (0, \e_\kappa )$. Then, the estimates
(\ref{eq:6.23}) and (\ref{eq:6.24}) in the above Lemma are enough to
show that
\[
\varphi \longmapsto \tilde{\mathcal N} (\e, h,k, \nu  ; \varphi )
\]
is a contraction from
\[ \{ \varphi  \in \tilde {\mathcal F} \quad : \quad \|
\varphi \|_{\tilde {\mathcal F}} \leq 2 \, c \, R_\e^{4-2m-\delta}
\} ,
\]
into itself and hence has a unique fixed point $\tilde \varphi_{\e,
h,k, \nu}$ in this set. This fixed point is a solution of
(\ref{eq:6.99}) in $N_{j,R_\e}$ and hence provides a constant scalar
curvature \K\ form on $N_{j, R_\e}$.

\medskip

We have obtained the~:
\begin{prop} There exist $c >0$ (independent of $\kappa$) and  $\e_\kappa =
\e (\kappa) > 0$ such that, for all $\e \in (0, \e_\kappa)$, for all
$h \in {\mathcal C}^{4, \alpha}(\del B_{\Gamma_j})$ and $k \in
{\mathcal C}^{2, \alpha}(\del B_{\Gamma_j})$ and $\nu \in{\mathbb
R}$ satisfying (\ref{eq:6.18}), the \K\ form
\[
\eta_{h , k  , \nu } : = \eta_j + i \, \del \, \bar \del \,
\varphi_{h , k  , \nu} ,
\]
defined on $N_{j ,R_\e}$, has constant scalar curvature equal to
$\e^{2} \, \nu$. Moreover
\[
\| \varphi_{h , k  , \nu} |_{C_{j ,R_\e/2} - C_{j , R_\e}} (R_\e \,
\cdot ) - H^{i}_{h , k } \|_{{\mathcal C}^{4, \alpha} (B_{j ,1}
-B_{j , 1/2})} \leq c \, R_\e^{ 4-2m} ,
\]
for some constant $c >0$ independent of $\kappa$ and $\nu$.
\label{pr:6.2}
\end{prop}
The important fact is that the last estimate involves a constant
times $R_\e^{4-2m}$  where the constant does not depend on $\kappa$
provided $\e \in (0, \e_\kappa)$.

\medskip

Using (\ref{eq:6.24}) and (\ref{eq:6.25}) and increasing $\tilde
c_\kappa$ if necessary, one checks that
\begin{equation}
\begin{array}{cllll} \| ( \tilde \varphi_{h,k, \nu }  - \tilde
\varphi_{h', k' , \nu'})\, |_{\bar C_{j , R_\e/2} -C_{j  , R_\e}}
(R_\e \, \cdot) - H_{h - h', k - k'}^{i}
\|_{{\mathcal C}^{4, \alpha} ( \bar B_{j ,1} - B_{j  ,1/2} )}  \\[3mm]
\hspace{60mm}  \leq \tilde c_\kappa \, ( R_\e^{\delta -1} \, \| (h -
h', k - k' )\|_{{\mathcal C}^{4, \alpha} \times {\mathcal C}^{2,
\alpha} } + R_\e^{4-2m} \, | \nu - \nu'| ).
\end{array}
\label{eq:6.29}
\end{equation}

\subsection{Cauchy data matching : the proof of Theorem~\ref{th:1.3}}

Building on the analysis of the previous sections we complete the
proof of Theorem~\ref{th:1.3}.

\medskip

Granted the results of Proposition~\ref{pr:6.1} and
Proposition~\ref{pr:6.2}, it remains to explain how to choose \[
{\bf h} : = (h_1, \ldots, h_n), \qquad {\bf k} : = (k_0, \ldots,
k_n)
\]
satisfying (\ref{eq:6.5}) and
\[ {\bf \tilde h} : = (\tilde h_1 ,
\ldots, \tilde h_n),  \qquad {\bf \tilde k} : = (\tilde k_1 ,
\ldots, \tilde k_n)
\]
satisfying (\ref{eq:6.18}) in such a way that, for each $j =1,
\ldots, n$, the function
\[
\psi_{j}^{o}  : = \left( \varphi_{j} + \varphi_{{\bf h},{\bf k}}
\right) (r_\e \, \cdot ),
\]
defined in $\bar B_{j  , 2 } - B_{j  , 1}$ (see
Proposition~\ref{pr:6.1}) on the one hand, and for
\[
\nu : = {\bf s} (\omega_{\e, {\bf h}, {\bf k}})
\]
the function
\[
\psi_j^{i}  : =  \e^{2} \, \left( \tilde \varphi_{j} + \tilde
\varphi_{\tilde h_j, \tilde k_j , {\bf s} (\omega_{h,k})} \right)
(R_\e \, \cdot ),
\]
defined in $\bar B_{j ,1} - B_{j  ,1/2}$ (see
Proposition~\ref{pr:6.1}) on the other hand, have their partial
derivatives up to order $3$ which coincide on $\del B_{j ,1}$.

\medskip

In any case, our aim is now to solve the following system of
equations
\begin{equation}
\psi_j^{o} =\psi_j^{i} , \qquad \del_r \, \psi_j^{o} = \del_r \,
\psi_j^{i} , \qquad \Delta_0 \, \psi_j^{o} = \Delta_0 \, \psi_j^{i},
\qquad \del_r \, \Delta_0  \, \psi_j^{o} = \del_r \, \Delta_0 \,
\psi_j^{i}, \label{eq:6.300}
\end{equation}
on $\del B_{j  ,1}$ where $r =|v|$ and $v= (v^{1}, \ldots, v^{m})$
are coordinates in $B_{j,2}$.

\medskip

Let us assume that we have already solved this problem. The first
identity in (\ref{eq:6.300}) implies that $\psi_j^{o}$ and
$\psi_j^{i}$ as well as all their $k$-th order partial derivatives
with respect any vector field tangent to $\del B_{j,1}$, with $k
\leq 4$, agree on $\del B_{j,1}$. The second identity in
(\ref{eq:6.300}) then shows that $\del_r \psi_j^{o}$ and $\del_r
\psi_j^{i}$ as well as all their $k$-th order partial derivatives
with respect any vector field tangent to $\del B_{j,1}$, with $k
\leq 3$, agree on $\del B_{j,1}$. Using the decomposition of the
Laplacian in polar coordinates, it is easy to check that the third
identity implies that $\del_r^2 \psi_j^{o}$ and $\del_r^2
\psi_j^{i}$ as well as all their $k$-th order partial derivatives
with respect any vector field tangent to $\del B_{j,1}$, with $k
\leq 2$, agree on $\del B_{j,1}$. And finally, the last identity in
(\ref{eq:6.300}) implies that $\del_r^3 \psi_j^{o}$ and $\del_r^3
\psi_j^{i}$ as well as all their first order partial derivative with
respect any vector field tangent to $\del B_{j,1}$, agree on $\del
B_{j,1}$.

\medskip

Moreover, the \K\ form
\[
i \, \del \, \bar \del \, ( \mbox{$\frac{1}{2}$} \, |v|^{2} +
\psi^{o}_j ),
\]
defined in $B_{j  ,2} - B_{j  ,1}$ and the \K\ form
\[
i \, \del \, \bar \del \, (\mbox{$\frac{1}{2}$} \, |v|^{2} +
\psi^{i}_j),
\]
defined in $B_{j, 1} - B_{j  , 1/2}$, both have the same constant
scalar curvature equal to ${\bf s} (\omega_{\e, {\bf h}, {\bf k}})$.
This then implies that any $k$-th order partial derivatives of the
functions $\psi^{o}_j$ and $\psi^{i}_j$, with $k \leq 4$, coincide
on $\del B_{j,1}$.

\medskip

Therefore, we conclude that the function $\psi$ defined by $\psi : =
\psi_j^{o}$ in $B_{j , 2} -B_{j , 1}$ and $\psi : = \psi_j^{i}$ in
$B_{j , 1} - B_{j , 1/2}$ is ${\mathcal C}^{4}$ in $B_{j,2} -
B_{j,1/2}$ and is a solution of the nonlinear elliptic partial
differential equation
\[
{\bf s} \,\left( i \, \del\, \bar \del (\mbox{$\frac{1}{2}$} \,
|v|^{2} + \psi ) \right) = {\bf s} (\omega_{\e, {\bf h}, {\bf k}}) =
cte .
\]
It then follows from elliptic regularity theory together with a
bootstrap argument that the function $\psi$ is in fact smooth.
Hence, by gluing the \K\ metrics $\omega_{{\bf h},{\bf k}}$ and
$\omega_{\tilde h_j, \tilde k_j}$ on the different pieces
constituting $M_{r_\e}$, we have produced a \K\ metric on $M_{r_\e}$
which has constant scalar curvature. This will end the proof of the
Theorem~\ref{th:1.3}.

\medskip

\begin{remar}
In dimension $2$, a slight modification is due since the functions
involve some $\log$ terms. In view of (\ref{eq:3.4bis}) and
(\ref{eq:6.07}), we consider the function $\tilde \psi^{i}_j$
defined by
\[ \psi_j^{i}  : =  \e^{2} \, \left( \tilde \varphi_{j} +
\tilde \varphi_{\tilde h_j, \tilde k_j , {\bf s} (\omega_{h,k})}
\right) (R_\e \, \cdot ) - \e^{2} \,  a_j \, \log R_\e +
\frac{k^{(0)}}{2} \, \log \, r_\e
\]
There is no loss of generality in doing so since changing the
potential by some constant function does not alter the corresponding
\K\ forms.
\end{remar}

\medskip

It remains to explain how to find the boundary data
\[
{\bf h} =(h_1, \ldots, h_n), \quad {\bf k} = (k_1, \ldots, k_n),
\quad {\bf \tilde h} = (\tilde h_1, \ldots, \tilde h_n) \qquad
\mbox{and} \qquad {\bf \tilde k} = (\tilde k_1, \ldots, \tilde k_n)
\]

We will make use of the following result~:
\begin{lemma}
Assume that $\Gamma$ is a discrete subgroup of $U(m)$ acting freely
on ${\mathbb C}^{m} -\{0\}$. The mapping
\[
\begin{array}{rclclll}
\mathcal P :& \mathcal C^{4,\alpha}(\del B_{\Gamma}) \times \mathcal
C^{2,\alpha}(\del B_{\Gamma}) & \longrightarrow & \mathcal
C^{3,\alpha}(\del B_{\Gamma} ) \times \mathcal C^{1,\alpha}(\del
B_{\Gamma} )
\\[3mm]
& (h,k)  &\longmapsto    & (\partial_{r} \, (H^{i}_{h, k}- H^{o}_{h,
k}), \partial_{r} \, \Delta_0 \, (H^{i}_{h, k}-  H^{o}_{h, k})) ,
\end{array}
\]
is an isomorphism. \label{le:6.3}
\end{lemma}
{\bf Proof~:} There are many different ways to prove this result
\cite{Fak-Pac} for example. Let us concentrate on the case where $m
\geq 3$ since the  case $m=2$ is essentially the same. We use the
formulas (\ref{eq:5.44}) and (\ref{eq:5.444}) to compute
\[
\partial_{r} \, (H^{i}_{h, k}-  H^{o}_{h, k}) = \sum_{\gamma=0}^{\infty} 2 \, (\gamma + m-1) \, \left( h^{(\gamma)}
+
\frac{k^{(\gamma)}}{(\gamma +m)(\gamma+m-2)} \, \right) \, e_\gamma
\]
and
\[
\partial_{r} \, \Delta_0 (H^{i}_{h, k}-  H^{o}_{h, k}) =
\sum_{\gamma=0}^{\infty} 2 \, (\gamma + m-1) \, k^{(\gamma)} \,
e_\gamma
\]
It is then easy to see that
\[
\begin{array}{rclclll}
{\mathcal P} :& W^{4,2}(\del B_{\Gamma}) \times W^{2,2}(\del
B_{\Gamma}) & \longrightarrow & W^{3,2}(\del B_{\Gamma} ) \times
W^{1,2}(\del B_{\Gamma} )
\\[3mm]
& (h,k)  &\longmapsto    & (\partial_{r} \, (H^{i}_{h, k}- H^{o}_{h,
k}), \partial_{r} \, \Delta_0 \, (H^{i}_{h, k}-  H^{o}_{h, k})) ,
\end{array}
\]
is well defined and invertible. Recall that the norm in $W^{\ell,
2}(\del B_\Gamma)$ can be taken to be \[ \| f \|_{W^{\ell, 2}} =
\left( \sum_{\gamma =0}^\infty ( 1+ \gamma)^{2\ell} \,
|f^{(\gamma)}|^2 \right)^{1/2}
\]
whenever the function $f$ is decomposed as
\[
f =  \sum_{\gamma =1}^\infty \, f^{(\gamma)} \, e_\gamma .
\]
Elliptic regularity theory then implies that the same result is true
when the operator is defined between H\"older spaces. \hfill $\Box$

\medskip

It will be convenient to observe that $\psi_{j}^{o}$ satisfies
\begin{equation}
\| \psi_{j}^{o} -  H_{h_j, k_j}\|_{{\mathcal C}^{4, \alpha}(B_{j ,
2} -B_{j ,1})}  \leq c \, r_\e^{4} , \label{eq:6.30}
\end{equation}
and also that
\begin{equation}
\| \psi_{j}^{i} -  \e^{2} \, \tilde H_{\tilde h_j, \tilde k_j}
\|_{{\mathcal C}^{4, \alpha}(B_{j , 1} - B_{j , 1/2})} \leq c \,
\e^{2} \, R_\e^{4-2n} = c \, r_\e^{4} , \label{eq:6.31}
\end{equation}
for some constant $c >0$ which does not depend on $\kappa$, provided
$\e$ is chosen small enough, say $\e \in (0, \e_\kappa)$. These two
estimates follow at once from the estimates in
Proposition~\ref{pr:6.1}, Proposition~\ref{pr:6.2} and also from the
choice of $r_\e$.

\medskip

We use the following notations for the rescaled boundary data
\[
({\bf h'}, {\bf k'} , {\bf \tilde h'} , {\bf \tilde k'}) : = ( {\bf
h}, {\bf k}, \e^{2} \, {\bf \tilde h} , \e^{2} \, {\bf \tilde k} ) .
\]
Using Lemma~\ref{le:6.3}, the solvability of (\ref{eq:6.300})
reduces to a fixed point problem which can be written as
\[
({\bf h'} ,  {\bf \tilde h'} , {\bf  k'}, {\bf \tilde k'} ) = S_\e (
{\bf h'} , {\bf \tilde h'} , {\bf k'} , {\bf \tilde k} ) ,
\]
and we know from (\ref{eq:6.30}) and (\ref{eq:6.31}) that the
nonlinear operator $S_\e$ satisfies
\[
\| S_\e ({\bf h'} ,  {\bf \tilde h'} , {\bf  k'}, {\bf \tilde k'}  )
\|_{({\mathcal C}^{4, \alpha})^{n} \times ({\mathcal C}^{2, \alpha}
)^{n}} \leq c_0 \, r_\e^{4} ,
\]
for some constant $c_0 >0$ which does not depend on $\kappa$,
provided $\e \in (0, \e_\kappa)$. We finally choose
\[
\kappa = 2 \, c_0 ,
\]
and $\e \in (0, \e_{\kappa})$. We have therefore proved that $S_\e$
is a map from
\[
A_\e : = \left\{ ({\bf h'} ,  {\bf \tilde h'} , {\bf  k'}, {\bf
\tilde k'}  ) \in ({\mathcal C}^{4, \alpha})^{n} \times ({\mathcal
C}^{2, \alpha})^{n} \quad : \quad \| ({\bf h'} ,  {\bf \tilde h'} ,
{\bf  k'}, {\bf \tilde k'} ) \|_{({\mathcal C}^{4, \alpha})^{n}
\times ({\mathcal C}^{2, \alpha})^{n} } \leq \kappa \, r_\e^{4}
\right\} ,
\]
into itself. It follows from (\ref{eq:6.14}), (\ref{eq:6.15}) and
(\ref{eq:6.29}) that, reducing $\e_\kappa$ if this is necessary,
$S_\e$ is a contraction mapping from $A_\e$ into itself for all $\e
\in (0, \e_\kappa)$. Therefore, $S_\e$ has a fixed point in this
set. This completes the proof of the existence of a solution of
(\ref{eq:6.300}).

\medskip

The proof of the existence on $M_{r_\e} $ of a \K\ form $\omega_\e$
which has constant scalar curvature is therefore complete. Observe
that the scalar curvature of $\omega$ and $\omega_\e$ are close
since the estimate
\[
|{\bf s} (\omega_\e) -{\bf s} (\omega)|\leq c \, r_\e^{2m}
\]
follows directly from the construction.

\section{Refined asymptotics for ALE spaces}

Let us now describe in detail $(N, \eta)$, the blow up at the origin
of ${\C}^{m}$ endowed with the Burns-Calabi-Simanca metric. Away
from the exceptional divisor, the \K \ form $\eta$ is given by
\[
\eta  = i \, \del \, \bar \del A_m (|v|^{2})
\]
where $v = (v_1, \ldots, v_n)$ are complex coordinates in ${\mathbb
C}^{m} -\{ 0 \}$ and where the function $s \longmapsto A_m(s)$ is a
solution of the ordinary differential equation
\[
s^{2} \, ( s \, \del_s A_m )^{m-1} \, \del_s^{2} A_m + (m-1) \,   s
\, \del_s A_m  -  (m-2) =0
\]
which satisfies $A_m \sim \log s$ near $0$. We refer to \cite{si}
for a derivation of this equation. It turns out that, when $m=2$,
the function $A_2$ is explicitly given by
\[
A_2 (s) =  \log  s + \lambda \, s
\]
where $\lambda >0$, while in dimension $m \geq 3$, even though there
is no explicit formula for $A_m$ we have the simple~:
\begin{lemma}
Assume that $m \geq 3$. Then the function $A_m$ can be expanded as
\[
A_m (s) = \lambda \, s - \lambda^{2-m} \, s^{2-m} + {\mathcal O}
(s^{1-m})
\]
for $s >1$, where $\lambda >0$. \label{le:7.1}
\end{lemma}
{\bf Proof~:} Define the function $\zeta$ by $ s \, \zeta : = s \,
\del_s A_m - 1$. A direct computation shows that $\zeta$ solves
\[
(1 + s \, \zeta)^{m-1} \, s^{2} \, \del_s \zeta =  (1+ s\,
\zeta)^{m-1} - 1 - (m-1) \, s \, \zeta
\]
If in addition we take $\zeta (0) = 1$, then $\del_s \zeta $ remains
positive and one can check that $\zeta$ is well defined for all time
and converges to some positive constant $\lambda$, as $s$ tends to
$\infty$. This immediately implies that $s \, \del_s A_m = \lambda
\, s + {\mathcal O} (1)$ at infinity. The expansion then follows
easily. \hfill $\Box$

\medskip

Changing variables $u : =  \sqrt{2  \, \lambda} \, v$, we see from
the previous Lemma that the \K\ form $\eta$ can be expanded near
infinity as
\begin{equation}
\eta = i \, \del \, \bar \del ( \mbox{$\frac{1}{2}$} \, |u|^{2} +
\log |u|^{2} ) , \label{eq:7.1}
\end{equation}
in dimension $m =2$ and as
\begin{equation}
\eta = i \, \del \, \bar \del ( \mbox{$\frac{1}{2}$} \, |u|^{2} -
2^{m-2} \, |u|^{4-2m} + {\mathcal O} (|u|^{3-2m})) , \label{eq:7.2}
\end{equation}
in dimension $m \geq 3$.

\medskip

We now recall the following results of Joyce \cite{j} (which is a
Corollary of his Theorem 8.2.3 in our notations):
\begin{teor}
\label{th:7.1} Given $\Gamma$ a finite subgroup of $SU(m)$ acting
freely on $\C^{m} - \{0\}$ and $\pi\colon X \rightarrow
\C^{m}/\Gamma$ a \K\ crepant resolution of $\C^{m}/\Gamma$. Then
there exists a Ricci-flat \K\ metric $\omega$ such that
\[
(\pi^{-1})^{*}\omega = i \, \del \, \bar \del ( \mbox{$\frac{1}{2}$}
\, |u|^{2} + \tilde{\varphi}(u))
\]
outside a compact neighborhood of $\pi^{-1}(0)$.

\medskip

Moreover
\[
\tilde{\varphi}(u) = |u|^{2-2m} + {\mathcal O} (|u|^{\gamma}))
\]
for some $\gamma \in (1-2m,2-2m)$.
\end{teor}

We end this section by a proof of Remark~\ref{re:3.1}.
\begin{lemma} Assume we are given a potential $\varphi$ defined on
$C_{\Gamma}$ such that $\varphi \in  {\mathcal C}^{4, \alpha}_{2
-\gamma} (\bar C_\Gamma)$, for some $\gamma >0$. Further assume that
\begin{equation}
\eta : = i \, \del \, \bar \del \, ( \mbox{$\frac{1}{2}$} \, |u|^2 +
\varphi ) \label{eq:7.3}
\end{equation}
is a zero scalar curvature \K\ form. Then, the function $ \varphi$
can be expanded as
\begin{equation}
\varphi = a \cdot u  + b + c \, |u|^{4-2m} + {\mathcal O}
(|u|^{3-2m}) , \label{eq:7.10}
\end{equation} when $m \geq 3$ and as
\begin{equation}
\varphi = a \cdot u + b + c \, \log |u| + {\mathcal O} (|u|^{-1}),
\label{eq:7.11} \end{equation} when $m =2$. Here $a\in {\mathbb C}$
and $b \in {\mathbb R}$. In particular, the potential  $\tilde
\varphi : = \varphi - a \cdot u - b$ satisfies
\[
\eta : = i \, \del \, \bar \del \, ( \mbox{$\frac{1}{2}$} \, |u|^{2}
+ \tilde \varphi ).
\]
\label{le:7.2}
\end{lemma}
{\bf Proof~:} The key point is that, since $\eta$ has zero scalar
curvature, the potential $\varphi$ is a solution of some nonlinear
fourth order elliptic differential equation and satisfies some {\it
a priori} bound. It is then possible to get "refined asymptotics"
for the potential $\varphi$ in the spirit of what has been done in
\cite{Kor-Maz-Pac-Sch} for constant scalar curvature metrics. These
refined asymptotics are obtained by using a bootstrap argument in
H\"older weighted spaces.

\medskip

Using (\ref{eq:4.0}), we see that  the scalar curvature of $\eta$
can be expanded in powers of $\varphi$ as
\[
{\bf s} \, (\omega) =  \frac{1}{2} \, \Delta_0^{2} \, \varphi +
Q_{g_{eucl}} (\nabla^{2} \varphi ),
\]
where the nonlinear operator $Q_{g_{eucl}}$ collects all the
nonlinear terms. We shall now be more specific about the structure
of $Q_{g_{eucl}}$. Indeed, it follows from the explicit computation
of the Ricci curvature that the nonlinear operator $Q_{g_{eucl}}$
can be decomposed as
\begin{equation}
Q_{g_{eucl}} (\nabla^{2} \, \varphi) =  \sum_q B_{q,4,2} (\nabla^{4}
\varphi , \nabla^{2} \varphi) \, C_{q, 4,2} (\nabla^{2} \varphi) +
\sum_q  B_{q,3,3} (\nabla^{3} \varphi , \nabla^{3} \varphi) \, C_{q,
3,3} (\nabla^{2} \varphi) ,\label{eq:7.4}
\end{equation}
where the sum over $q$ is finite, the operators $(U,V)
\longrightarrow B_{q,a,b} (U,V)$ are bilinear in the entries and
have coefficients which are bounded functions in ${\mathcal C}^{0,
\alpha} (\bar C_{\Gamma})$. The nonlinear operators $W
\longrightarrow C_{q,a,b} (W)$ have Taylor expansion (with respect
to $W$) whose coefficients are bounded functions on ${\mathcal
C}^{0, \alpha} (\bar C_{\Gamma})$.

\medskip

If we assume that $\varphi \in  {\mathcal C}^{4, \alpha}_{2 -\gamma}
(\bar C_\Gamma)$, then we see that
\[
Q_{g_{eucl}} (\nabla^{2} \varphi) \in  {\mathcal C}^{0, \alpha}_{-2
-2 \gamma} (\bar C_\Gamma)
\]
Therefore, $\Delta_0^{2} \, \varphi \in {\mathcal C}^{0, \alpha}_{-2
-2 \gamma} (\bar C_\Gamma)$.

\medskip

Now, if $\Delta^{2}_0 \, \varphi \in {\mathcal C}^{0,
\alpha}_{\gamma'-4}  (\bar C_\Gamma)$ and $\varphi \in  {\mathcal
C}^{4, \alpha}_{2 -\gamma} (\bar C_\Gamma)$ for some $\gamma
>0$ then, depending on the value of $\gamma'$, the following
alternative hold \cite{Maz}~:
\begin{itemize}

\item[(i)] If $\gamma' \in (1,2)$, then $\varphi \in {\mathcal
C}^{4, \alpha}_{\gamma'} (\bar C_\Gamma)$.\\

\item[(ii)] If $\gamma' \in (0,1)$, then  $\varphi \in {\mathcal
C}^{4, \alpha}_{\gamma'} (\bar C_\Gamma) \oplus \{ u \longmapsto
a\cdot u \quad : \quad  a \in {\mathbb C}\}$.\\

\item[(iii)] If $m \geq 3$ and $\gamma' \in (4-2m ,0)$, then
\[
\varphi \in {\mathcal C}^{4, \alpha}_{\gamma'} (\bar C_\Gamma)
\oplus \{ u \longmapsto a\cdot u \quad  : \quad a \in {\mathbb C}\}
\oplus {\mathbb R}.
\]

\item[(iv)] If $m \geq 3$ and $\gamma' \in (3-2m,4-2m)$, then
\[
\varphi \in {\mathcal C}^{4, \alpha}_{\gamma'} (\bar C_\Gamma)
\oplus \{ u \longmapsto a\cdot u \quad : \quad  a \in {\mathbb C}\}
\oplus {\mathbb R} \oplus \mbox{Span}\{ u \longmapsto |u|^{4-2m} \}.
\]

\item[(v)] If $m = 2$ and $\gamma' \in (-1,0)$, then
\[
\varphi \in {\mathcal C}^{4, \alpha}_{\gamma'} (\bar C_\Gamma)
\oplus \{ u \longmapsto a\cdot u \quad : \quad  a \in {\mathbb C}\}
\oplus {\mathbb R} \oplus \mbox{Span}\{ u \longmapsto \log |u| \}.
\]

\end{itemize}
Using these together with a bootstrap argument, we conclude that
(\ref{eq:7.10}) and (\ref{eq:7.11}) hold. The result then follows by
taking $\tilde \varphi :=  \varphi  - a \cdot u -b$. \hfill $\Box$

\section{Applications, examples and comments}

{\bf Blow up of smooth manifolds :} Theorem~\ref{th:1.1} follows at
once from Theorem~\ref{th:1.3} and the analysis of
Lemma~\ref{le:7.1} by taking $(N_j, \eta_j)= (N, a_j \, \eta)$ where
$(N,\eta)$ is the Blow up at the origine of ${\mathbb C}^m$ endowed
with the Burns-Calabi-Simanca metric and $a_j >0$. Observe that the
points of blow up $p_1, \ldots, p_n$ and the coefficients $a_1,
\ldots, a_n$ are parameters of our construction.

\medskip

A first natural question is to which base smooth manifolds can
Theorem~\ref{th:1.1} be applied ! Here, we do not make a
comprehensive list but we highlight some large class of manifolds~:
\begin{enumerate}
\item[(i)]
All the \K -Einstein manifold with discrete automorphism group. This
means any manifold with negative first Chern class and many families
of examples of positive first Chern class we know of \cite{t1},
\cite{n}, \cite{agp}. We should note that there are no \K -Einstein
manifolds {\em Futaki nondegenerate} except the ones with discrete
automorphisms as observed by LeBrun-Simanca \cite{ls}. \\

\item[(ii)]
Most of the zero scalar curvature \K\ surfaces which have been
proved by Kim, LeBrun Pontecorvo, Rollin and Singer \cite{lbsing},
\cite{lb2}, \cite{klbp}, \cite{kp} to admit such constant scalar
curvature metric. In particular any blow up of a non Ricci-flat \K\
surface whose integral of the scalar curvature is non-negative has
blow ups which admit zero scalar curvature \K\ metrics. Of course if
the number of blow ups is sufficiently large no continuous families
of
automorphisms survive and we can then apply Theorem~\ref{th:1.3}.\\

\item[(iii)]
Note that also flat tori of any dimension can be used as base
manifolds, since, despite the presence of continuous automorphisms,
there are no nonzero holomorphic vector fields vanishing somewhere.
Their first Chern class being zero, Corollary~\ref{co:1.2} does not
apply.\\

\item[(iv)]
Some important classes of manifolds on which there are constant
scalar curvature \K\ metrics have been provided by Fine \cite{f}.
Indeed, he has proved existence of \K\ constant scalar curvature
metrics on complex surfaces with a holomorphic submersion onto a
Riemann surface $\Sigma$ with smooth fibres of genus at least two.
If the genus of $\Sigma $ is larger than or equal to $2$, the
automorphism group is indeed discrete.\\

\item[(v)]
Another family of examples of  constant scalar curvature \K\
manifolds with discrete automorphism group has been given by Hong
\cite{h1}, \cite{h2}. These are ruled manifolds given by the
projectivization of some vector bundles over constant scalar
curvature \K\ manifolds. \\

\item[(vi)]
In \cite{ls}, LeBrun and Simanca gave examples (and strategies to
construct new ones) of {\em Futaki nondegenerate} manifolds with
constant scalar curvature \K\ metrics.\\

\item[(vii)]
Recall that the space of holomorphic vector fields on a blow
manifold is isomorphic to the space of those holomorphic vector
fields on the base manifold vanishing at the blow up points. Hence
our procedure applied to any of the {\em nondegenerate} manifold
above gives new {\em nondegenerate} manifolds (with constant scalar
curvature
by our result), so our procedure can be iterated.\\

\item[(viii)]
Riemannian products of {\em nondegenerate} \K\  manifolds of
constant scalar curvature is again a {\em nondegenerate} \K\
manifold of constant scalar curvature. By taking factors with scalar
curvature of different signs and scaling one can then produce also
on the blow ups \K\ metrics of any nonzero scalar curvature.

\end{enumerate}

In addition, to any of the above examples, one can apply
LeBrun-Simanca's implicit function argument \cite{ls} to get open
subset of the \K\ cone of fixed complex manifolds and also open
subset of moduli of variations of complex structures for which
constant scalar curvature \K\ metrics exist, providing a wealth of
new examples.

\medskip

{\bf Zero scalar curvature examples. Proof of Corollary
\ref{co:1.1}~:} Let us now focus on the effect of our construction
on the size of the scalar curvature when we blow up smooth points.
Let us denote by $\pi$ the standard projection from the blow up
manifold $\tilde M$ to the base manifold $M$. To this aim let us
recall that the average of the scalar curvature of a \K\ metric is a
cohomological number given  by
\[
{\bf s} (\omega) = \frac{m c_1(M) \cup
[\omega]^{m-1}([M])}{[\omega]^{m} ([M])} .
\]
Our gluing procedure constructs on $\tilde M$ metrics in the \K\
classes
\[
[\omega_{\e}] = \pi^{*}[\omega] - \e^{2} (a_1 \, PD[E_1]+ \ldots +
a_n \, PD[E_n]) \] while the first Chern class behaves like
\[
c_1(\tilde M) = \pi^{*}(c_1(M)) - (m-1) \, (PD[E_1] + \ldots +PD
[E_n]) .\]

Recalling (see \cite{Gri-Har} page 475) that for any $j=1, \dots,
n$,
\[ (PD[E_{j}])^{m}[\tilde{M}] = (-1)^{m-1}, \] we get

\[
(c_1(\tilde M))\cup
([\omega_{\e}])^{m-1} ([\tilde{M}])  = \\
(c_1(M)\cup[\omega]^{m-1})[M] - \e^{2m-2}(m-1)(\sum_{j=1}^{n}a_j)
\]

and
\[
[\omega_{\e}]^{m}([\tilde{M}]) = [\omega]^{m}([M]) +
(-1)^{m-1}\e^{2m}(\sum_{j=1}^{n} a_j)
\]

The scalar curvature of this metric is hence given by
\[
{\bf s}(\omega_{\e}) = m \frac{(c_1(M)\cup[\omega]^{m-1})([M]) -
\e^{2m-2}(m-1)(\sum_{j=1}^{n}a_j)} {[\omega]^{m}([M]) +
(-1)^{m-1}\e^{2m}(\sum_{j=1}^{n} a_j)} .
\]
It is easily seen that, since $a_{j}>0$, this gives a decreasing
function of $\e$,  for $\e$ close to $0$ (and of course it gives the
old scalar curvature for $\e = 0$).

\medskip

The direct application of Theorem \ref{th:1.1} would then give small
negative scalar curvature if $(M,\omega)$ had zero scalar curvature.
Nonetheless changing the \K\ class we can bypass this problem
provided the first Chern class of the base orbifold is nonzero,
prescribing the scalar curvature to vanish in the gluing procedure.

\begin{corol}
Any blow up (at a finite set of smooth points) of a compact smooth
\K\ manifold (or orbifold) of zero scalar curvature of discrete type
with nonzero first Chern class, has a \K\ metric of zero constant
scalar curvature.
\end{corol}
{\bf Proof~:} Let us denote by $\omega(0)$ the  zero scalar
curvature \K\ metric on the base manifold $M$, and by $\rho$ the
harmonic representative of the first Chern class $c_1(M)$ (hence non
zero by our assumption). LeBrun-Simanca have proved  (\cite{lbsm},
Corollary $1$) that, if the first Chern class is nonzero, the
automorphism group is discrete and for $|t|$ is sufficiently small
(say $t\in [-t_0,t_0]$), each \K\ class $[ \omega (0) - t \rho]$
contains a metric $\omega (t)$ of constant scalar curvature and this
constant is positive for $t > 0$ and negative for $t < 0$. Moreover
$\omega (t)$ depends continuously on $t$.

\medskip

We can apply Theorem~\ref{th:1.3} to the continuous family of \K\
forms $\omega (t)$. Given $t \in [-t_0, t_0]$, this yields the
existence of $\e_0(t) >0$ and a family of \K\ metrics $\omega (t
,\e)$ of constant scalar curvature for all $\e  \in (0,\e_0(t))$. It
turns out that, the constant $\e_0(t)$ are uniformly bounded from
below by some positive constant $\e_0 >0$ since $\e_0(t)$ only
depends on the ${\mathcal C}^{2, \alpha}$ norm of the coefficients
of the \K\ form $\omega (t)$ and these are uniformly bounded as $t
\in [-t_0, t_0]$. We claim that, reducing $\e_0$ if this is
necessary, $\omega(t, \e)$ depends continuously on $t$. This follows
easily from the fact that the \K\ forms on the blown up manifold are
obtained by solving nonlinear problems using a fixed point theorems
for a contraction mapping. Therefore, they depend continuously on
any of the parameters of our construction such as the \K\ class, the
parameter $\e$, the points which are blown up, the coefficients $a_j
>0$, \ldots

\medskip

Let us then look at the family of constant scalar curvature metrics
$\omega (t , \e)$. We known that, reducing $\e_0$ if necessary,
$\omega (-t_0 ,\e )$ has constant negative scalar curvature while
$\omega (t_0 , \e )$ has positive scalar curvature, for all $\e \in
(0, \e_0)$. Moreover ${\bf s} (\omega (t , \e))$ depends
continuously on $t$. Therefore, for each $\e\in (0, \e_0)$, there
exists $t_\e \in [-t_0,t_0]$ such that ${\bf s}(\omega (t_\e , \e ))
=0$  as claimed. \hfill $\Box$

\medskip

\noindent Note that the above Corollary can be applied to most of
the examples described above in (ii) and (viii).

\medskip

{\bf Desingularization of orbifolds~:} More delicate is the
situation for singular orbifolds since few examples even of \K
-Einstein orbifolds are known. As mentioned in the introduction, the
clearest picture is in the complex dimension $2$ and $3$, where,
thanks to the work of Kronheimer \cite{kr} and Joyce \cite{j} we
know how to handle $SU(m)$ singular points. We summarize this in the
following

\begin{corol}
\label{222} Let $(M, \omega)$ be a {\em nondegenerate} compact
$m$-dimensional constant scalar curvature \K\ orbifold with $m=2$ or
$3$ and isolated singularities. Let $p_1, \dots p_n \in M$ be any
set of points with a neighborhood biholomorphic to a neighborhood of
the origin in ${\mathbb C}^{m} / \, \Gamma_j$, where $\Gamma_j$ is a
finite subgroup of $SU(m)$. Let further $N_j$ be a \K\ crepant
resolution of ${\mathbb C}^{m} / \, \Gamma_j$ (which always exists,
see \cite{bpvdv} for $m=2$ and \cite{ro} for $m=3$), and $\eta_j$
given by Theorem \ref{th:7.1}.

\medskip

Then there exists $\e_0 > 0$, such that, for all $\e \in (0, \e_0)$,
there exists a constant scalar curvature \K\ form $\omega_\e$ on $M
\sqcup _{{p_{1}, \e}} N_1 \sqcup_{{p_{2}, \e}} \dots \sqcup
_{{p_{n}, \e}} N_n$.

\medskip

Moreover,
\begin{itemize}
\item[(i)]
if $\omega$ had positive (resp. negative) scalar curvature then
$\omega_{\e}$ has positive (resp. negative) scalar curvature, \\
\item[(ii)]
if $c_1(M)\neq 0$ and $\omega_M$ has zero scalar curvature then
$\omega_{\e}$ can be chosen to have zero scalar curvature too.
\end{itemize}
\end{corol}

\medskip

The range of applicability of Corollary \ref{222} is very large,
even if we look just at \K -Einstein orbifolds of non positive
scalar curvature thanks to Aubin-Yau's solution of the Calabi
conjecture (which holds in the orbifold category). In fact we can
use it to prove the following general result mentioned in the
introduction

\begin{corol}
Any compact complex surface of general type admits constant scalar
curvature \K\ metrics. \label{co:8.3}
\end{corol}

The proof of the above result requires some notions from algebraic
geometry which can be found for example in \cite{bpvdv} and which we
quickly recall for reader's convenience.

\medskip

First of all a complex surface $M$ is called {\em minimal} if it
does not contain  a smooth rational curve of self-intersection $-1$.
A fundamental result in complex surface theory (the
Enriques-Castelnuovo Criterion, see e.g. \cite{Gri-Har} page 476)
says that any such curve is in fact the exceptional divisor of a
blow up at a smooth point of a smooth surface. Moreover one can
apply the above procedure (``blowing down") a finite number of times
to be left with a minimal surface uniquely defined, hence called
{the minimal model} $\bar{M}$ of $M$.

\medskip

From a different perspective one can study an algebraic surface by
looking at its images into projective spaces, via maps given by
evaluating holomorphic sections of line bundles as in the celebrated
Kodaira's embedding theorem. In particular, if $K_{M}$ is the
canonical line bundle of $M$, one has (rational) maps
$\phi_{K_M^{\otimes k}}$ from $M$ into $\mathbb{P}(H^{0}(M,
K_M^{\otimes k}))$. These in general may not be defined at points
which annihilates all holomorphic sections of $K_M^{\otimes k}$, but
for minimal surfaces of general type they are indeed globally defined
holomorphic maps for $k\geq 5$ (see e.g. \cite{bpvdv}, page 220).

\medskip

A complex surface $M$ is said to be {\em of general type} if $dim
(\phi_{K_M^{\otimes k}}(M)) = 2$ for $k$ large enough. If $M$ is a minimal 
surface of general type
Kodaira \cite{kod} proved that $\phi_{K_M^{\otimes k}}$ is an
embedding away from smooth rational curves of self-intersection
$-2$, and Brieskorn \cite{br} has proved that the image of these
curves are isolated singular points of the image surface with local
structure groups $\Gamma_j$, with $\Gamma_j \subset SU(2)$. We are
now in position to give the proof of the above corollary.

\medskip

{\bf Proof of Corollary~\ref{co:8.3}~:} Let us first assume $M$ is a
minimal complex surface of general type and suppose $k$ is chosen
big enough to guarantee that the image of the the pluricanonical
rational map $\phi_{K_M^{\otimes k}}$ is an embedding away from the set of
$(-2)$-curves of $M$, which get collapsed to points, giving the
singularities of $\phi_{K_M^{\otimes k}}$.

\medskip

Kobayashi \cite{k} has proved that $\phi_{K_M^{\otimes k}}$ has a \K
-Einstein orbifold metric of negative scalar curvature, extending
Aubin's proof of the Calabi conjecture. Moreover $c_1(M) < 0$
implies, as in the smooth case, the existence of only a discrete
group of automorphisms.
\medskip

As already observed, being the structure groups of the singularities
in $SU(2)$, we have an ALE local model with the required decay at
infinity. We can then apply Theorem~\ref{th:1.3}. The complex
manifold produced by our gluing construction is easily seen to be
minimal, hence
getting a constant negative scalar curvature \K\ metric on the
minimal resolution $Y$ of $M$. But $M$ is already a minimal model of
$M$, therefore the minimal model of $M$, and so $M$ is in fact $Y$
proving our result.

\medskip

If $M$ is not minimal, we apply the previous discussion to its
minimal model $Y$, which is a complex surface with discrete
automorphism group to get a \K\ constant negative scalar curvature
metric. Recalling that $M$ is obtained from $Y$ applying a finite
number of blow ups, Theorem~\ref{th:1.1} (possibly applied more than
once in case one needs to blow up at a point on the exceptional
divisor of the previous blow up, and of course blowing up preserves
the property of having only discrete automorphism groups) gives the
conclusion. \hfill $\Box$

\medskip

Going back to the problem of resolving singularities in the \K\
constant scalar curvature setting, in dimension greater than $3$
only few examples can be dealt at the moment.

\medskip

Other types of singularities which can be dealt with are, for
example, those locally modeled on ${\mathbb C}^{m} /\intero_m$,
where $\intero_m$ acts diagonally on ${\mathbb C}^{m}$, by
multiplication by a fixed the $m$-th root of unity $\zeta =
e^{\frac{2\pi i}{m}}$. Putting $r = (|z^{1}|^{2} + \cdots +
|z^{m}|^{2})^{1/2}$, Calabi \cite{ca} defined a \K\ potential on $X
\setminus \{exceptional \,\, divisor\}$ by
\[
\varphi = (r^{2m} + 1)^{\frac{1}{m}} + {\frac{1}{m}}
\sum_{j=0}^{m-1} \zeta^{j} \log ((r^{2m} + 1)^{\frac{1}{m}}  -
\zeta^{j})\,\, .
\]
We can then observe that $\eta = \frac{i}{2} \, \partial
 \, {\bar{\partial}} \, \varphi$ is indeed a \K\ form which extends
through the exceptional divisor and is ALE, Ricci flat, asymptotic
to ${\mathbb C}^{m}/\intero_m$. We can then glue $(X, \eta)$ to any
smooth \K\ orbifold $(M,\omega)$ of constant scalar curvature
provided the Futaki obstructions vanishes as described in Section 4.

\medskip

The above example has been recently generalized by Rollin-Singer
\cite{rs2}. They have shown that if $G=\{1, \lambda, \dots,
\lambda^{k-1}\}$, $\lambda = e^{\frac{2 \pi i}{k}}$, then ${\mathbb
C}^{m}/G$ has an ALE scalar flat (in general not Ricci-flat) \K\
resolution whose metric decays at infinity of order $2-2m$.

These last examples can be used to produce compact orbifolds
by taking global quotients of some of the smooth manifolds described
in the first section (e.g. tori or \K -Einstein manifolds with
negative first Chern class or with positive first Chern class and
discrete automorphism group containing a group as above).

\end{document}